\documentclass[12pt,vatola]{article}
\usepackage{amsfonts}

\usepackage{mathrsfs}

\usepackage{graphicx}

\def\c{\centerline}

\def\re#1{\par\hangindent\parindent\indent\llap{#1\enspace}\ignorespaces}

\def\no{\noindent}

\begin{document}

\c{\bf\large Geometrical Theory on Combinatorial Manifolds}

\vskip 5mm

\c{Linfan Mao}\vskip 2mm

\c{\scriptsize (Chinese Academy of Mathematics and System Science,
Beijing 100080, P.R.China)}

\c{\scriptsize E-mail: maolinfan@163.com}

\vskip 6mm

\begin{minipage}{130mm}

\no{\bf Abstract}: {\small For an integer $m\geq 1$, a combinatorial
manifold $\widetilde{M}$ is defined to be a geometrical object
$\widetilde{M}$ such that for $\forall p\in\widetilde{M}$, there is
a local chart $(U_p,\varphi_p)$ enable $\varphi_p:U_p\rightarrow
B^{n_{i_1}}\bigcup B^{n_{i_2}}\bigcup\cdots\bigcup B^{n_{i_{s(p)}}}$
with $B^{n_{i_1}}\bigcap B^{n_{i_2}}\bigcap\cdots\bigcap
B^{n_{i_{s(p)}}}\not=\emptyset$, where $B^{n_{i_j}}$ is an
$n_{i_j}$-ball for integers $1\leq j\leq s(p)\leq m$. Topological
and differential structures such as those of $d$-pathwise connected,
homotopy classes, fundamental $d$-groups in topology and tangent
vector fields, tensor fields, connections, {\it Minkowski} norms in
differential geometry on these finitely combinatorial manifolds are
introduced. Some classical results are generalized to finitely
combinatorial manifolds. Euler-Poincare characteristic is discussed
and geometrical inclusions in Smarandache geometries for various
geometries are also presented by the geometrical theory on finitely
combinatorial manifolds in this paper. }\vskip 2mm

\no{\bf Key Words}: {\small manifold, finitely combinatorial
manifold, topological structure, differential structure,
combinatorially Riemannian geometry, combinatorially Finsler
geometry, Euler-Poincare characteristic.}\vskip 2mm

\no{\bf AMS(2000)}: {\small 51M15, 53B15, 53B40, 57N16}

\end{minipage}

\vskip 5mm

\no{\bf \S $1.$ Introduction}

\vskip 3mm

\no As a model of spacetimes in physics, various geometries such as
those of Euclid, Riemannian and Finsler geometries are established
by mathematicians. Today, more and more evidences have shown that
our spacetime is not homogenous. Thereby models established on
classical geometries are only unilateral. Then {\it are there some
kinds of overall geometries for spacetimes in physics?} The answer
is YES. Those are just Smarandache geometries established in last
century but attract more one's attention now. According to the
summary in $[4]$, they are formally defined following.

\vskip 3mm

\no{\bf Definition $1.1$}($[4],[17]$) \ {\it A Smarandache geometry
is a geometry which has at least one Smarandachely denied
axiom($1969$), i.e., an axiom behaves in at least two different ways
within the same space, i.e., validated and invalided, or only
invalided but in multiple distinct ways.

A Smarandache $n$-manifold is a $n$-manifold that support a
Smarandache geometry.}

\vskip 3mm

For verifying the existence of Smarandache geometries, Kuciuk and
Antholy gave a popular and easily understanding example on an Euclid
plane in $[4]$. In $[3]$, Iseri firstly presented a systematic
construction for Smarandache geometries by equilateral triangular
disks on Euclid planes, which are really Smarandache $2$-dimensional
geometries (see also $[5]$). In references $[6],[7]$ and $[13]$,
particularly in $[7]$, a general constructing way for Smarandache
$2$-dimensional geometries on maps on surfaces, called {\it map
geometries} was introduced, which generalized the construction of
Iseri. For the case of dimensional number$\geq 3$, these {\it
pseudo-manifold geometries} are proposed, which are approved to be
Smarandache geometries and containing these Finsler and K\"{a}hler
geometries as sub-geometries in $[12]$.

In fact, by the Definition $1.1$ a general but more natural way for
constructing Smarandache geometries should be seeking for them on a
union set of spaces with an axiom validated in one space but
invalided in another, or invalided in a space in one way and another
space in a different way. These unions are so called Smarandache
multi-spaces. This is the motivation for this paper. Notice that in
$[8]$, these {\it multi-metric spaces} have been introduced, which
enables us to constructing Smarandache geometries on multi-metric
spaces, particularly, on multi-metric spaces with a same metric.

\vskip 3mm

\no{\bf Definition $1.2$} \ {\it A multi-metric space
$\widetilde{A}$ is a union of spaces $A_1,A_2,\cdots,A_m$ for an
integer $k\geq 2$ such that each $A_i$ is a space with metric
$\rho_i$ for $\forall i, 1\leq i\leq m$.}

\vskip 3mm

Now for any integer $n$, these $n$-manifolds $M^n$ are the main
objects in modern geometry and mechanics, which are locally
euclidean spaces ${\bf R}^n$ satisfying the $T_2$ separation axiom
in fact, i.e., for $\forall p,q\in M^n$, there are local charts
$(U_p,\varphi_p)$ and $(U_q,\varphi_q)$ such that $U_p\bigcap
U_q=\emptyset$ and $\varphi_p: U_p\rightarrow {\bf B}^n$,
$\varphi_q: U_q\rightarrow {\bf B}^n$, where

$$B^n=\{(x_1,x_2,\cdots,x_n)| x_1^2+x_2^2+\cdots+x_n^2 \ < 1\}.$$

\no is an open ball.

These manifolds are locally euclidean spaces. In fact, they are also
homogenous spaces. But the world is not homogenous. Whence, a more
important thing is considering these combinations of different
dimensions, i.e., {\it combinatorial manifolds} defined following
and finding their good behaviors for mathematical sciences besides
just to research these manifolds. Two examples for these
combinations of manifolds with different dimensions in ${\bf R}^3$
are shown in Fig.$1.1$, in where, (a) represents a combination of a
$3$-manifold, a torus and $1$-manifold, and (b) a torus with $4$
bouquets of $1$-manifolds.\vskip 2mm

\includegraphics[bb=25 10 200 115]{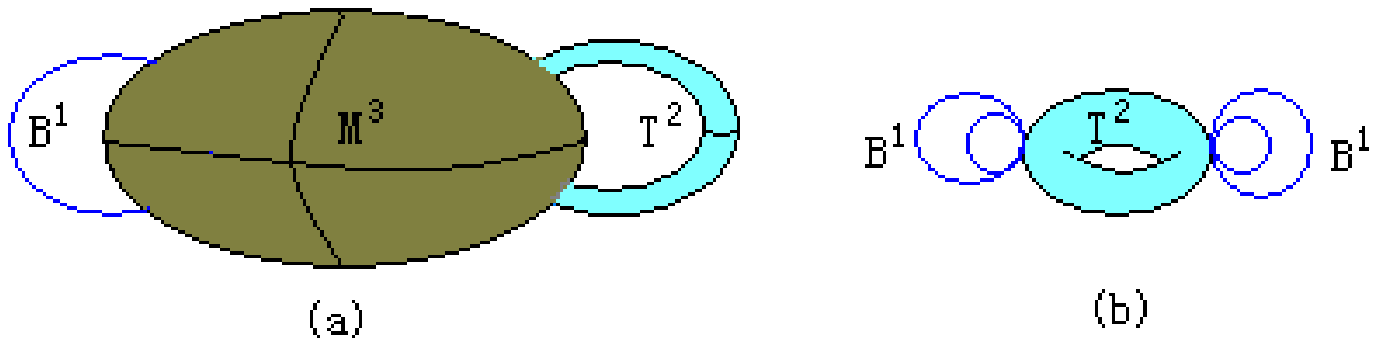}

\c{\bf Fig.$1.1$}\vskip 2mm

For an integer $s\geq 1$, let $n_1,n_2,\cdots,n_s$ be an integer
sequence with $0< n_1< n_2<\cdots< n_s$. Choose $s$ open unit balls
$B_1^{n_1}, B_2^{n_2}, \cdots,B_s^{n_s}$, where
$\bigcap\limits_{i=1}^s B_i^{n_i}\not=\emptyset$ in ${\bf
R}^{n_1+_2+\cdots n_s}$. Then a {\it unit open combinatorial ball of
degree $s$} is a union

$$\widetilde{B}(n_1, n_2,\cdots,n_s) = \bigcup\limits_{i=1}^sB_i^{n_i}.$$

\vskip 3mm

\no{\bf Definition $1.3$}\ {\it For a given integer sequence
$n_1,n_2,\cdots,n_m, m\geq 1$ with $0< n_1< n_2<\cdots< n_s$, a
combinatorial manifold $\widetilde{M}$ is a Hausdorff space such
that for any point $p\in \widetilde{M}$, there is a local chart
$(U_p,\varphi_p)$ of $p$, i.e., an open neighborhood $U_p$ of $p$ in
$\widetilde{M}$ and a homoeomorphism $\varphi_p:
U_p\rightarrow\widetilde{B}(n_1(p), n_2(p),\cdots,n_{s(p)}(p))$ with
$\{n_1(p),
n_2(p),\cdots,n_{s(p)}(p)\}\subseteq\{n_1,n_2,\cdots,n_m\}$ and
$\bigcup\limits_{p\in\widetilde{M}}\{n_1(p),
n_2(p),\cdots,n_{s(p)}(p)\}=\{n_1,n_2,\cdots,n_m\}$, denoted by
$\widetilde{M}(n_1,n_2,\cdots,n_m)$ or $\widetilde{M}$ on the
context and

$$\widetilde{{\mathcal A}}=\{(U_p,\varphi_p)|
p\in\widetilde{M}(n_1,n_2,\cdots,n_m))\}$$

\no an atlas on $\widetilde{M}(n_1,n_2,\cdots,n_m)$. The maximum
value of $s(p)$ and the dimension $\widehat{s}(p)$ of
$\bigcap\limits_{i=1}^{s(p)}B_i^{n_i}$ are called the dimension and
the intersectional dimensional of $\widetilde{M}(n_1,n_2,$
$\cdots,n_m)$ at the point $p$, respectively.

A combinatorial manifold $\widetilde{M}$ is called finite if it is
just combined by finite manifolds. }

Notice that $\bigcap\limits_{i=1}^s B_i^{n_i}\not=\emptyset$ by the
definition of unit combinatorial balls of degree $s$. Thereby, for
$\forall p\in\widetilde{M}(n_1,n_2,\cdots,n_s)$, either it has a
neighborhood $U_p$ with $\varphi_p:U_p\rightarrow{\bf
R}^{\varsigma}$, $\varsigma\in\{n_1,n_2,\cdots,n_s\}$ or a
combinatorial ball $\widetilde{B}(\tau_1,\tau_2,\cdots,\tau_l)$ with
$\varphi_p:U_p\rightarrow\widetilde{B}(\tau_1,\tau_2,\cdots,\tau_l)$,
$l\leq s$ and
$\{\tau_1,\tau_2,\cdots,\tau_l\}\subseteq\{n_1,n_2,\cdots,n_s\}$
hold.

The main purpose of this paper is to characterize these finitely
combinatorial manifolds, such as those of topological behaviors and
differential structures on them by a combinatorial method. For these
objectives, topological and differential structures such as those of
$d$-pathwise connected, homotopy classes, fundamental $d$-groups in
topology and tangent vector fields, tensor fields, connections, {\it
Minkowski} norms in differential geometry on these combinatorial
manifolds are introduced. Some results in classical differential
geometry are generalized to finitely combinatorial manifolds. As an
important invariant,  Euler-Poincare characteristic is discussed and
geometrical inclusions in Smarandache geometries for various
existent geometries are also presented by the geometrical theory on
finitely combinatorial manifolds in this paper.

For terminologies and notations not mentioned in this section, we
follow $[1]-[2]$ for differential geometry, $[5],[7]$ for graphs and
$[14],[18]$ for topology.

\vskip 6mm

\no{\bf \S $2.$ Topological structures on combinatorial manifolds}

\vskip 4mm

\no By a topological view, we introduce topological structures and
characterize these finitely combinatorial manifolds in this section.

\vskip 4mm

\no{\bf $2.1.$ Pathwise connectedness}

\vskip 3mm

\no On the first, we define {\it $d$-dimensional pathwise
connectedness} in a finitely combinatorial manifold for an integer
$d, d\geq 1$, which is a natural generalization of {\it pathwise
connectedness} in a topological space.

\vskip 4mm

\no{\bf Definition $2.1$} \ {\it For two points $p,q$ in a finitely
combinatorial manifold $\widetilde{M}(n_1,n_2,$ $\cdots,n_m)$, if
there is a sequence $B_1,B_2,\cdots,B_s$ of $d$-dimensional open
balls with two conditions following hold.}\vskip 3mm

(1) \ {\it $B_i\subset\widetilde{M}(n_1,n_2,\cdots,n_m)$ for any
integer $i, 1\leq i\leq s$ and $p\in B_1$, $q\in B_s$;}

(2) \ {\it The dimensional number ${\rm dim}(B_i\bigcap B_{i+1})\geq
d$ for $\forall i, 1\leq i\leq s-1$.}\vskip 2mm

\no{\it Then points $p, q$ are called $d$-dimensional connected in
$\widetilde{M}(n_1,n_2,\cdots,n_m)$ and the sequence
$B_1,B_2,\cdots,B_e$ a $d$-dimensional path connecting $p$ and
$q$, denoted by $P^d(p,q)$.

If each pair $p,q$ of points in the finitely combinatorial manifold
$\widetilde{M}(n_1,n_2,\cdots,n_m)$ is $d$-dimensional connected,
then $\widetilde{M}(n_1,n_2,\cdots,n_m)$ is called $d$-pathwise
connected and say its connectivity$\geq d$.}

\vskip 3mm

Not loss of generality, we consider only finitely combinatorial
manifolds with a connectivity$\geq 1$ in this paper. Let
$\widetilde{M}(n_1,n_2,\cdots,n_m)$ be a finitely combinatorial
manifold and $d, d\geq 1$ an integer. We construct a labelled graph
$G^d[\widetilde{M}(n_1,n_2,\cdots,n_m)]$ by

$$V(G^d[\widetilde{M}(n_1,n_2,\cdots,n_m)])= V_1\bigcup V_2,$$

\no where $V_1=\{n_i-{\rm manifolds} \ M^{n_i} \ {\rm in} \
\widetilde{M}(n_1,n_2,\cdots,n_m) |1\leq i\leq m\}$ and $V_2=\{ {\rm
isolated \ intersection \ points} \ O_{M^{n_i},M^{n_j}} \ {\rm of}
M^{n_i}, M^{n_j} \ {\rm in} \ \widetilde{M}(n_1,n_2,\cdots,n_m) \
{\rm for} \ 1\leq i,j\leq m\}$. Label $n_i$ for each $n_i$-manifold
in $V_1$ and $0$ for each vertex in $V_2$ and

$$E(G^d[\widetilde{M}(n_1,n_2,\cdots,n_m)])= E_1\bigcup E_2,$$

\no where \ $E_1=\{(M^{n_i},M^{n_j})| {\rm dim}(M^{n_i}\bigcap
M^{n_j})\geq d, 1\leq i,j\leq m\}$   \ and
 \ $E_2=\{(O_{M^{n_i},M^{n_j}},M^{n_i}), (O_{M^{n_i},M^{n_j}},M^{n_j})|
M^{n_i} \ {\rm tangent} \ M^{n_j} \ {\rm at \ the \ point} \
O_{M^{n_i},M^{n_j}}$ $ {\rm for} \ 1\leq i,j\leq m\}$.

\includegraphics[bb=25 10 200 180]{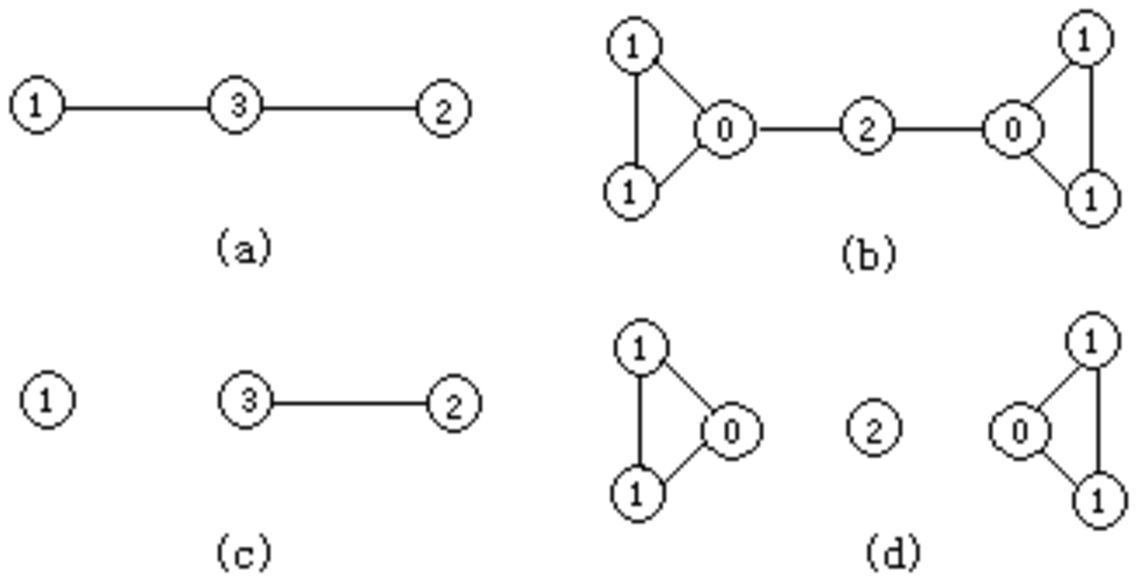}\vskip 2mm

\c{\bf Fig.$2.1$}\vskip 2mm

For example, these correspondent labelled graphs gotten from
finitely combinatorial manifolds in Fig.$1.1$ are shown in
Fig.$2.1$, in where $d=1$ for (a) and (b), $d=2$ for (c) and (d). By
this construction, properties following can be easily gotten.

\vskip 4mm

\no{\bf Theorem $2.1$} \ {\it Let
$G^d[\widetilde{M}(n_1,n_2,\cdots,n_m)]$ be a labelled graph of a
finitely combinatorial manifold $\widetilde{M}(n_1,n_2,\cdots,n_m)$.
Then}

($1$) {\it $G^d[\widetilde{M}(n_1,n_2,\cdots,n_m)]$ is connected
only if $d\leq n_1$.}

($2$) {\it there exists an integer $d, d\leq n_1$ such that
$G^d[\widetilde{M}(n_1,n_2,\cdots,n_m)]$ is connected.}

\vskip 3mm

{\it Proof} \ By definition, there is an edge $(M^{n_i},M^{n_j})$
in $G^d[\widetilde{M}(n_1,n_2,\cdots,n_m)]$ for $1\leq i,j\leq m$
if and only if there is a $d$-dimensional path $P^d(p,q)$
connecting two points $p\in M^{n_i}$ and $q\in M^{n_j}$. Notice
that

$$(P^d(p,q)\setminus M^{n_i})\subseteq M^{n_j} \ {\rm and} \
(P^d(p,q)\setminus M^{n_j})\subseteq M^{n_i}.$$

\no Whence,

$$d\leq\min\{n_i,n_j\}. \hskip 60mm (2.1)$$

Now if $G^d[\widetilde{M}(n_1,n_2,\cdots,n_m)]$ is connected, then
there is a $d$-path $P(M^{n_i}, M^{n_j})$ connecting vertices
$M^{n_i}$ and $M^{n_j}$ for $\forall M^{n_i}, M^{n_j}\in
V(G^d[\widetilde{M}(n_1,n_2,\cdots,n_m)])$. Not loss of generality,
assume

$$P(M^{n_i}, M^{n_j})=M^{n_i}M^{s_1}M^{s_2}\cdots M^{s_{t-1}}M^{n_j}.$$

\no Then we get that

$$d\leq\min\{n_i,s_1,s_2,\cdots, s_{t-1},n_j\} \hskip 36mm (2.2)$$

\no by $(2.1)$. However, according to Definition $1.4$ we know that

$$\bigcup\limits_{p\in\widetilde{M}}\{n_1(p),
n_2(p),\cdots,n_{s(p)}(p)\}=\{n_1,n_2,\cdots,n_m\}. \ \ (2.3)$$

\no Therefore, we get that

$$d\leq\min(\bigcup\limits_{p\in\widetilde{M}}\{n_1(p),
n_2(p),\cdots,n_{s(p)}(p)\})=\min\{n_1,n_2,\cdots,n_m\}=n_1$$

\no by combining $(2.3)$ with $(2.3)$. Notice that points labelled
with $0$ and $1$ are always connected by a path. We get the
conclusion ($1$).

For the conclusion ($2$), notice that any finitely combinatorial
manifold is always pathwise $1$-connected by definition.
Accordingly, $G^1[\widetilde{M}(n_1,n_2,\cdots,n_m)]$ is connected.
Thereby, there at least one integer, for instance $d=1$ enabling
$G^d[\widetilde{M}(n_1,n_2,$ $\cdots,n_m)]$ to be connected. This
completes the proof. \ \ \ $\natural$

According to Theorem $2.1$, we get immediately two corollaries
following.

\vskip 4mm

\no{\bf Corollary $2.1$} \ {\it For a given finitely combinatorial
manifold $\widetilde{M}$, all connected graphs $G^d[\widetilde{M}]$
are isomorphic if $d\leq n_1$, denoted by $G[\widetilde{M}]$.}

\vskip 3mm

\no{\bf Corollary $2.2$} \ {\it If there are $k$ $1$-manifolds
intersect at one point $p$ in a finitely combinatorial manifold
$\widetilde{M}$, then there is an induced subgraph $K^{k+1}$ in
$G[\widetilde{M}]$.}

\vskip 3mm

Now we define an edge set $E^d(\widetilde{M})$ in $G[\widetilde{M}]$
by

$$E^d(\widetilde{M})=
E(G^d[\widetilde{M}])\setminus E(G^{d+1}[\widetilde{M}]).$$

\no Then we get a graphical recursion equation for graphs of a
finitely combinatorial manifold $\widetilde{M}$ as a by-product.

\vskip 4mm

\no{\bf Theorem $2.2$}\ {\it Let $\widetilde{M}$ be a finitely
combinatorial manifold. Then for any integer $d, d\geq 1$, there is
a recursion equation

$$G^{d+1}[\widetilde{M}]=G^d[\widetilde{M}]
-E^d(\widetilde{M})$$

\no for graphs of $\widetilde{M}$.}

\vskip 3mm

{\it Proof} \ It can be obtained immediately by definition. \ \
$\natural$

For a given integer sequence $1\leq n_1< n_2<\cdots< n_m, m\geq 1$,
denote by ${\mathcal H}^d(n_1,n_2,\cdots,n_m)$ all these finitely
combinatorial manifolds $\widetilde{M}(n_1,n_2,\cdots,n_m)$ with
connectivity$\geq d$, where $d\leq n_1$ and ${\mathcal
G}(n_1,n_2,\cdots,n_m)$ all these connected graphs
$G[n_1,n_2,\cdots,n_m]$ with vertex labels $0,n_1,n_2,\cdots,n_m$
and conditions following hold.\vskip 3mm

($1$) \ The induced subgraph by vertices labelled with $1$ in $G$
is a union of complete graphs;

($2$) \ All vertices labelled with $0$ can only be adjacent to
vertices labelled with $1$.\vskip 2mm

Then we know a relation between sets ${\mathcal
H}^d(n_1,n_2,\cdots,n_m)$ and ${\mathcal G}(n_1,n_2,\cdots,n_m)$.

\vskip 4mm

\no{\bf Theorem $2.3$} \ {\it Let $1\leq n_1<n_2<\cdots<n_m, m\geq
1$ be a given integer sequence. Then every finitely combinatorial
manifold $\widetilde{M}\in{\mathcal H}^d(n_1,n_2,\cdots,n_m)$
defines a labelled connected graph
$G[n_1,n_2,\cdots,n_m]\in{\mathcal G}(n_1,n_2,\cdots,n_m)$.
Conversely, every labelled connected graph
$G[n_1,n_2,\cdots,n_m]\in{\mathcal G}(n_1,n_2,\cdots,n_m)$ defines a
finitely combinatorial manifold $\widetilde{M}\in{\mathcal
H}^d(n_1,n_2,$ $\cdots,n_m)$ for any integer $1\leq d\leq n_1$.}

\vskip 3mm

{\it Proof} \ For $\forall\widetilde{M}\in{\mathcal
H}^d(n_1,n_2,\cdots,n_m)$, there is a labelled graph
$G[n_1,n_2,\cdots,n_m]\in{\mathcal G}(n_1,n_2,\cdots,n_m)$
correspondent to $\widetilde{M}$ is already verified by Theorem
$2.1$. For completing the proof, we only need to construct a
finitely combinatorial manifold $\widetilde{M}\in{\mathcal
H}^d(n_1,n_2,\cdots,n_m)$ for $\forall
G[n_1,n_2,\cdots,n_m]\in{\mathcal G}(n_1,n_2,\cdots,n_m)$. Denoted
by $l(u)=s$ if the label of a vertex $u\in V(G[n_1,n_2,\cdots,n_m])$
is $s$. The construction is carried out by the following
programming.\vskip 2mm

\no STEP $1.$ \ Choose $|G[n_1,n_2,\cdots,n_m]|-|V_0|$ manifolds
correspondent to each vertex $u$ with a dimensional $n_i$ if
$l(u)=n_i$, where $V_0=\{u|u\in V(G[n_1,n_2,\cdots,n_m]) \ {\rm and}
\ l(u)=0\}$. Denoted by $V_{\geq 1}$ all these vertices in
$G[n_1,n_2,\cdots,n_m]$ with label$\geq 1$.\vskip 2mm

\no STEP $2.$ \ For $\forall u_1\in V_{\geq 1}$ with
$l(u_1)=n_{i_1}$, if its neighborhood set
$N_{G[n_1,n_2,\cdots,n_m]}(u_1)\bigcap$ $ V_{\geq
1}=\{v_1^1,v_1^2,\cdots,v_1^{s(u_1)}\}$ with $l(v_1^1)=n_{11}$,
$l(v_1^2)=n_{12}$, $\cdots$, $l(v_1^{s(u_1)})=n_{1s(u_1)}$, then let
the manifold correspondent to the vertex $u_1$ with an intersection
dimension$\geq d$ with manifolds correspondent to vertices
$v_1^1,v_1^2,\cdots,v_1^{s(u_1)}$ and define a vertex set
$\Delta_1=\{u_1\}$.\vskip 2mm

\no STEP $3.$ \ If the vertex set
$\Delta_l=\{u_1,u_2,\cdots,u_l\}\subseteq V_{\geq 1}$ has been
defined and $V_{\geq 1}\setminus\Delta_l\not=\emptyset$, let
$u_{l+1}\in V_{\geq 1}\setminus\Delta_l$ with a label $n_{i_{l+1}}$.
Assume

$$(N_{G[n_1,n_2,\cdots,n_m]}(u_{l+1})\bigcap V_{\geq 1})
\setminus\Delta_l=\{v_{l+1}^1,v_{l+1}^2,\cdots,v_{l+1}^{s(u_{l+1})}\}$$

\no with $l(v_{l+1}^1)=n_{l+1,1}$, $l(v_{l+1}^2)=n_{l+1,2}$,
$\cdots$,$l(v_{l+1}^{s(u_{l+1})})=n_{l+1,s(u_{l+1})}$. Then let the
manifold correspondent to the vertex $u_{l+1}$ with an intersection
dimension$\geq d$ with manifolds correspondent to these vertices
$v_{l+1}^1,v_{l+1}^2, \cdots,v_{l+1}^{s(u_{l+1})}$ and define a
vertex set $\Delta_{l+1}=\Delta_l\bigcup\{u_{l+1}\}$.\vskip 2mm

\no STEP $4.$ \ Repeat steps $2$ and $3$ until a vertex set
$\Delta_{t}=V_{\geq 1}$ has been constructed. This construction is
ended if there are no vertices $w\in V(G)$ with $l(w)=0$, i.e.,
$V_{\geq 1}=V(G)$. Otherwise, go to the next step.\vskip 2mm

\no STEP $5.$ \ For $\forall w\in
V(G[n_1,n_2,\cdots,n_m])\setminus V_{\geq 1}$, assume
$N_{G[n_1,n_2,\cdots,n_m]}(w)=\{w_1,w_2,$ $\cdots,w_e\}$. Let all
these manifolds correspondent to vertices $w_1,w_2,\cdots,w_e$
intersects at one point simultaneously and define a vertex set
$\Delta^*_{t+1}=\Delta_t\bigcup\{w\}$.\vskip 2mm

\no STEP $6.$ \ Repeat STEP $5$ for vertices in
$V(G[n_1,n_2,\cdots,n_m])\setminus V_{\geq 1}$. This construction
is finally ended until a vertex set
$\Delta^*_{t+h}=V(G[n_1,n_2,\cdots,n_m])$ has been
constructed.\vskip 2mm

As soon as the vertex set $\Delta^*_{t+h}$ has been constructed, we
get a finitely combinatorial manifold $\widetilde{M}$. It can be
easily verified that $\widetilde{M}\in{\mathcal
H}^d(n_1,n_2,\cdots,n_m)$ by our construction way. \quad\quad
$\natural$

\vskip 5mm

\no{\bf $2.2$ Combinatorial equivalence}

\vskip 3mm

\no For a finitely combinatorial manifold $\widetilde{M}$ in
${\mathcal H}^d(n_1,n_2,\cdots,n_m)$, denoted by
$G[\widetilde{M}(n_1,$ $n_2,\cdots,n_m)]$ and $G[\widetilde{M}]$ the
correspondent labelled graph in ${\mathcal G}(n_1,n_2,\cdots,n_m)$
and the graph deleted labels on
$G[\widetilde{M}(n_1,n_2,\cdots,n_m)]$, $C(n_i)$ all these vertices
with a label $n_i$ for $1\leq i\leq m$, respectively.

\vskip 4mm

\no{\bf Definition $2.2$} \ {\it Two finitely combinatorial
manifolds $\widetilde{M}_1(n_1,n_2,\cdots,n_m),
\widetilde{M}_2(k_1,k_2,$ $\cdots,k_l)$ are called equivalent if
these correspondent labelled graphs}

$$G[\widetilde{M}_1(n_1,n_2,\cdots,n_m)]\cong G[\widetilde{M}_2(k_1,k_2,\cdots,k_l)].$$

\vskip 3mm

Notice that if $\widetilde{M}_1(n_1,n_2,\cdots,n_m),
\widetilde{M}_2(k_1,k_2,\cdots,k_l)$ are equivalent, then we can get
that $\{n_1,n_2,\cdots,n_m\}=\{k_1,k_2,\cdots,k_l\}$ and
$G[\widetilde{M}_1]\cong G[\widetilde{M}_2]$. Reversing this idea
enables us classifying finitely combinatorial manifolds in
${\mathcal H}^d(n_1,n_2,\cdots,n_m)$ by the action of automorphism
groups of these correspondent graphs without labels.

\vskip 4mm

\no{\bf Definition $2.3$} \ {\it A labelled connected graph
$G[\widetilde{M}(n_1,n_2,\cdots,n_m)]$ is combinatorial unique if
all these correspondent finitely combinatorial manifolds
$\widetilde{M}(n_1,n_2,\cdots,n_m)$ are equivalent.}

\vskip 3mm

A labelled graph $G[n_1,n_2,\cdots,n_m]$ is called {\it
class-transitive} if the automorphism group ${\rm Aut}G$ is
transitive on $\{C(n_i), 1\leq i\leq m\}$. We find a characteristic
for combinatorially unique graphs.

\vskip 4mm

\no{\bf Theorem $2.4$} \ {\it A labelled connected graph
$G[n_1,n_2,\cdots,n_m]$ is combinatorially unique if and only if it
is class-transitive.}

\vskip 3mm

{\it Proof} \ For two integers $i,j, 1\leq i,j\leq m$, re-label
vertices in $C(n_i)$ by $n_j$ and vertices in $C(n_j)$ by $n_i$ in
$G[n_1,n_2,\cdots,n_m]$. Then we get a new labelled graph
$G'[n_1,n_2,\cdots,n_m]$ in ${\mathcal G}[n_1,n_2,\cdots,n_m]$.
According to Theorem $2.3$, we can get two finitely combinatorial
manifolds $\widetilde{M}_1(n_1,n_2,\cdots,n_m)$ and
$\widetilde{M}_2(k_1,k_2,\cdots,k_l)$ correspondent to
$G[n_1,n_2,\cdots,n_m]$ and $G'[n_1,n_2,\cdots,n_m]$.

Now if $G[n_1,n_2,\cdots,n_m]$ is combinatorially unique, we know
$\widetilde{M}_1(n_1,n_2,\cdots,n_m)$ is equivalent to
$\widetilde{M}_2(k_1,k_2,\cdots,k_l)$, i.e., there is an
automorphism $\theta\in {\rm Aut}G$ such that
$C^{\theta}(n_i)=C(n_j)$ for $\forall i,j, 1\leq i,j\leq m$.

On the other hand, if $G[n_1,n_2,\cdots,n_m]$ is class-transitive,
then for integers $i,j, 1\leq i,j\leq m$, there is an automorphism
$\tau\in{\rm Aut}G$ such that $C^{\tau}(n_i)=C(n_j)$. Whence, for
any re-labelled graph $G'[n_1,n_2,\cdots,n_m]$, we find that

$$G[n_1,n_2,\cdots,n_m]\cong G'[n_1,n_2,\cdots,n_m],$$

\no which implies that these finitely combinatorial manifolds
correspondent to $G[n_1,n_2,$ $\cdots,n_m]$ and
$G'[n_1,n_2,\cdots,n_m]$ are combinatorially equivalent, i.e.,
$G[n_1,n_2,\cdots,n_m]$ is combinatorially unique.  \ \ \ \
$\natural$

Now assume that for parameters $t_{i1},t_{i2},\cdots,t_{is_i}$, we
have known an enufunction

$$C_{M^{n_i}}[x_{i1},x_{i2},\cdots]
=\sum\limits_{t_{i1},t_{i2},\cdots,t_{is}}
n_i(t_{i1},t_{i2},\cdots,t_{is})x_{i1}^{t_{i1}}x_{i2}^{t_{i2}}\cdots
x_{is}^{t_{is}}$$

\no for $n_i$-manifolds, where $n_i(t_{i1},t_{i2},\cdots,t_{is})$
denotes the number of non-homeomorphic $n_i$-manifolds with
parameters $t_{i1},t_{i2},\cdots, t_{is}$. For instance the
enufunction for compact $2$-manifolds with parameter genera is

$$C_{\widetilde{M}}[x](2)= 1+\sum\limits_{p\geq 1}2x^p.$$

\no Consider the action of ${\rm Aut}G[n_1,n_2,\cdots,n_m]$ on
$G[n_1,n_2,\cdots,n_m]$. If the number of orbits of the automorphism
group ${\rm Aut}G[n_1,n_2,\cdots,n_m]$ action on $\{C(n_i), 1\leq
i\leq m\}$ is $\pi_0$, then we can only get $\pi_0!$ non-equivalent
combinatorial manifolds correspondent to the labelled graph
$G[n_1,n_2,\cdots,n_m]$ similar to Theorem $2.4$. Calculation shows
that there are $l!$ orbits action by its automorphism group for a
complete $(s_1+s_2+\cdots+s_l)$-partite graph
$K(k_1^{s_1},k_2^{s_2},\cdots,k_l^{s_l})$, where $k_i^{s_i}$ denotes
that there are $s_i$ partite sets of order $k_i$ in this graph for
any integer $i, 1\leq i\leq l$, particularly, for
$K(n_1,n_2,\cdots,n_m)$ with $n_i\not=n_j$ for $i,j, 1\leq i,j\leq
m$, the number of orbits action by its automorphism group is $m!$.
Summarizing all these discussions, we get an enufunction for these
finitely combinatorial manifolds $\widetilde{M}(n_1,n_2,\cdots,n_m)$
correspondent to a labelled graph $G[n_1,n_2,\cdots,n_m]$ in
${\mathcal G}(n_1,n_2,\cdots,n_m)$ with each label$\geq 1$.

\vskip 4mm

\no{\bf Theorem $2.5$} \ {\it Let $G[n_1,n_2,\cdots,n_m]$ be a
labelled graph in ${\mathcal G}(n_1,n_2,\cdots,n_m)$ with each
label$\geq 1$. For an integer $i, 1\leq i\leq m$, let the
enufunction of non-homeomorphic $n_i$-manifolds with given
parameters $t_1,t_2,\cdots,$ be $C_{M^{n_i}}[x_{i1},x_{i2},\cdots]$
and $\pi_0$ the number of orbits of the automorphism group ${\rm
Aut}G[n_1,n_2,\cdots,n_m]$ action on $\{C(n_i), 1\leq i\leq m\}$,
then the enufunction of combinatorial manifolds
$\widetilde{M}(n_1,n_2,\cdots,n_m)$ correspondent to a labelled
graph $G[n_1,n_2,\cdots,n_m]$ is}

$$C_{\widetilde{M}}(\overline{x})=
\pi_0!\prod\limits_{i=1}^mC_{M^{n_i}}[x_{i1},x_{i2},\cdots],$$

\no{\it particularly, if
$G[n_1,n_2,\cdots,n_m]=K(k_1^{s_1},k_2^{s_2},\cdots,k_m^{s_m})$
such that the number of partite sets labelled with $n_i$ is $s_i$
for any integer $i, 1\leq i\leq m$, then the enufunction
correspondent to $K(k_1^{s_1},k_2^{s_2},\cdots,k_m^{s_m})$ is}

$$C_{\widetilde{M}}(\overline{x})=m!
\prod\limits_{i=1}^mC_{M^{n_i}}[x_{i1},x_{i2},\cdots]$$

\no{\it and the enufunction correspondent to a complete graph
$K_m$ is}

$$C_{\widetilde{M}}(\overline{x})=
\prod\limits_{i=1}^mC_{M^{n_i}}[x_{i1},x_{i2},\cdots].$$

\vskip 3mm

{\it Proof} \ Notice that the number of non-equivalent finitely
combinatorial manifolds correspondent to $G[n_1,n_2,\cdots,n_m]$ is

$$\pi_0\prod\limits_{i=1}^mn_i(t_{i1},t_{i2},\cdots,t_{is})$$

\no for parameters $t_{i1},t_{i2},\cdots,t_{is}, 1\leq i\leq m$ by
the product principle of enumeration. Whence, the enufunction of
combinatorial manifolds $\widetilde{M}(n_1,n_2,\cdots,n_m)$
correspondent to a labelled graph $G[n_1,n_2,\cdots,n_m]$ is

\begin{eqnarray*}
C_{\widetilde{M}}(\overline{x})&=&
\sum\limits_{t_{i1},t_{i2},\cdots,t_{is}}
(\pi_0\prod\limits_{i=1}^mn_i(t_{i1},t_{i2},\cdots,t_{is}))
\prod\limits_{i=1}^mx_{i1}^{t_{i1}}x_{i2}^{t_{i2}}\cdots
x_{is}^{t_{is}}\\
&=& \pi_0!\prod\limits_{i=1}^mC_{M^{n_i}}[x_{i1},x_{i2},\cdots]. \ \
\ \natural
\end{eqnarray*}

\vskip 5mm

\no{\bf $2.3$ Homotopy classes}

\vskip 3mm

\no Denote by $f\simeq g$ two homotopic mappings $f$ and $g$.
Following the same pattern of homotopic spaces, we define
homotopically combinatorial manifolds in the next.

\vskip 4mm

\no{\bf Definition $2.4$} \ {\it Two finitely combinatorial
manifolds $\widetilde{M}(k_1,k_2,\cdots,k_l)$ and
$\widetilde{M}(n_1,n_2,$ $\cdots,n_m)$ are said to be homotopic if
there exist continuous maps

$f:\widetilde{M}(k_1,k_2,\cdots,k_l)\rightarrow\widetilde{M}(n_1,n_2,\cdots,n_m),$

$g:\widetilde{M}(n_1,n_2,\cdots,n_m)\rightarrow\widetilde{M}(k_1,k_2,\cdots,k_l)$

\no such that
$gf\simeq$identity$:\widetilde{M}(k_1,k_2,\cdots,k_l)
\rightarrow\widetilde{M}(k_1,k_2,\cdots,k_l)$ and
$fg\simeq$identity$:\widetilde{M}(n_1,n_2,\cdots,n_m)
\rightarrow\widetilde{M}(n_1,n_2,\cdots,n_m)$.}

\vskip 3mm

For equivalent homotopically combinatorial manifolds, we know the
following result under these correspondent manifolds being
homotopic. For this objective, we need an important lemma in
algebraic topology.

\vskip 4mm

\no{\bf Lemma $2.1$}(Gluing Lemma, [16]) \ {\it Assume that a space
$X$ is a finite union of closed subsets:
$X=\bigcup\limits_{i=1}^nX_i$. If for some space $Y$, there are
continuous maps $f_i: X_i\rightarrow Y$ that agree on overlaps,
i.e., $f_i|_{X_i\bigcap X_j}=f_j|_{X_i\bigcap X_j}$ for all $i,j$,
then there exists a unique continuous $f: X\rightarrow Y$ with
$f|_{X_i}=f_i$ for all $i$.}

\vskip 3mm

\no{\bf Theorem $2.6$} \ {\it Let
$\widetilde{M}(n_1,n_2,\cdots,n_m)$ and
$\widetilde{M}(k_1,k_2,\cdots,k_l)$ be finitely combinatorial
manifolds with an equivalence
$\varpi:G[\widetilde{M}(n_1,n_2,\cdots,n_m)]\rightarrow
G[\widetilde{M}(k_1,k_2,\cdots,k_l)]$. If for $\forall M_1,M_2\in
V(G[\widetilde{M}(n_1,n_2,\cdots,n_m)])$, $M_i$ is homotopic to
$\varpi(M_i)$ with homotopic mappings
$f_{M_i}:M_i\rightarrow\varpi(M_i)$, $g_{M_i}:\varpi(M_i)\rightarrow
M_i$ such that $f_{M_i}|_{M_i\bigcap M_j}=f_{M_j}|_{M_i\bigcap
M_j}$, $g_{M_i}|_{M_i\bigcap M_j}=g_{M_j}|_{M_i\bigcap M_j}$
providing $(M_i,M_j)\in E(G[\widetilde{M}(n_1,n_2,\cdots,n_m)])$ for
$1\leq i,j\leq m$, then $\widetilde{M}(n_1,n_2,\cdots,n_m)$ is
homotopic to $\widetilde{M}(k_1,k_2,\cdots,k_l)$.}

\vskip 3mm

{\it Proof} \ By the Gluing Lemma, there are continuous mappings

\c{$f:\widetilde{M}(n_1,n_2,\cdots,n_m)\rightarrow
\widetilde{M}(k_1,k_2,\cdots,k_l)$}

\no and

\c{$g:\widetilde{M}(k_1,k_2,\cdots,k_l)\rightarrow
\widetilde{M}(n_1,n_2,\cdots,n_m)$}

\no such that

\c{$f|_M=f_M \ {\rm and} \ g|_{\varpi(M)}=g_{\varpi(M)}$}

\no for $\forall M\in V(G[\widetilde{M}(n_1,n_2,\cdots,n_m)])$.
Thereby, we also get that

\c{$gf\simeq \ identity:\widetilde{M}(k_1,k_2,\cdots,k_l)
\rightarrow\widetilde{M}(k_1,k_2,\cdots,k_l)$}

\no and

\c{$fg\simeq \ identity:\widetilde{M}(n_1,n_2,\cdots,n_m)
\rightarrow\widetilde{M}(n_1,n_2,\cdots,n_m)$}

\no as a result of $g_Mf_M\simeq \ identity: M\rightarrow M$,
$f_Mg_M\simeq \ identity: \varpi(M)\rightarrow\varpi(M)$. \
$\natural$

We have known that a finitely combinatorial manifold
$\widetilde{M}(n_1,n_2,\cdots,n_m)$ is $d$-pathwise connected for
some integers $1\leq d\leq n_1$. This consequence enables us
considering fundamental $d$-groups of finitely combinatorial
manifolds.

\vskip 4mm

\no{\bf Definition $2.5$} \ {\ Let
$\widetilde{M}(n_1,n_2,\cdots,n_m)$ be a finitely combinatorial
manifold. For an integer $d, 1\leq d\leq n_1$ and $\forall
x\in\widetilde{M}(n_1,n_2,\cdots,n_m)$, a fundamental $d$-group at
the point $x$, denoted by
$\pi^d(\widetilde{M}(n_1,n_2,\cdots,n_m),x)$ is defined to be a
group generated by all homotopic classes of closed $d$-pathes based
at $x$.}

\vskip 3mm

If $d=1$ and $\widetilde{M}(n_1,n_2,\cdots,n_m)$ is just a manifold
$M$, we get that

$$\pi^d(\widetilde{M}(n_1,n_2,\cdots,n_m),x)=\pi(M,x).$$

\no Whence, fundamental $d$-groups are a generalization of
fundamental groups in topology. We obtain the following
characteristics for fundamental $d$-groups of finitely combinatorial
manifolds.

\vskip 4mm

\no{\bf Theorem $2.7$} \ {\it Let
$\widetilde{M}(n_1,n_2,\cdots,n_m)$ be a $d$-connected finitely
combinatorial manifold with $1\leq d\leq n_1$. Then}\vskip 3mm

($1$) \ {\it for $\forall x\in\widetilde{M}(n_1,n_2,\cdots,n_m)$,}

$$
\pi^d(\widetilde{M}(n_1,n_2,\cdots,n_m),x)\cong
(\bigoplus\limits_{M\in V(G^d)}\pi^d(M))\bigoplus \pi(G^d),
$$

\no{\it where $G^d=G^d[\widetilde{M}(n_1,n_2,\cdots,n_m)]$,
$\pi^d(M), \pi(G^d)$ denote the fundamental $d$-groups of a manifold
$M$ and the graph $G^d$, respectively and}

($2$) \ {\it for $\forall
x,y\in\widetilde{M}(n_1,n_2,\cdots,n_m)$,}

$$\pi^{d}(\widetilde{M}(n_1,n_2,\cdots,n_m),x)\cong\pi^{d}(\widetilde{M}(n_1,n_2,\cdots,n_m),y).$$

\vskip 3mm

{\it Proof} \ For proving the conclusion (1), we only need to prove
that for any cycle $\widetilde{C}$ in
$\widetilde{M}(n_1,n_2,\cdots,n_m)$, there are elements
$C_1^M,C_2^M,\cdots,C_{l(M)}^M\in\pi^d(M)$,
$\alpha_1,\alpha_2,\cdots,\alpha_{\beta(G^d)}\in\pi(G^d)$ and
integers $a_i^M, b_j$ for $\forall M\in V(G^d)$ and $1\leq i\leq
l(M)$, $1\leq j\leq c(G^d)\leq\beta(G^d)$ such that

$$\widetilde{C}\equiv\sum\limits_{M\in V(G^d)}
\sum\limits_{i=1}^{l(M)}a_i^MC_i^M+\sum\limits_{j=1}^{c(G^d)}b_j\alpha_j({\rm
mod} 2)$$

\no and it is unique. Let $C_1^M,C_2^M,\cdots,C_{b(M)}^M$ be a
base of $\pi^d(M)$ for $\forall M\in V(G^d)$. Since
$\widetilde{C}$ is a closed trail, there must exist integers
$k_i^M, l_j, 1\leq i\leq b(M), 1\leq j\leq\beta(G^d)$ and $h_P$
for an open $d$-path on $\widetilde{C}$ such that

$$\widetilde{C}=\sum\limits_{M\in V(G^d)}
\sum\limits_{i=1}^{b(M)}k_i^MC_i^M+\sum\limits_{j=1}^{\beta(G^d)}l_j\alpha_j+\sum\limits_{P\in
\Delta}h_PP,$$

\no where $h_P\equiv 0({\rm mod}2)$ and $\Delta$ denotes all of
these open $d$-paths on  $\widetilde{C}$. Now let

$$\{a_i^M|1\leq i\leq l(M)\}=\{k_i^M| k_i^M\not=0 \ {\rm and} \ 1\leq i\leq b(M) \},$$

$$\{b_j| 1\leq j\leq c(G^d)\}=\{l_j| l_j\not=0, 1\leq j\leq\beta(G^d) \}.$$

\no Then we get that

$$\widetilde{C}\equiv\sum\limits_{M\in V(G^d)}
\sum\limits_{i=1}^{l(M)}a_i^MC_i^M+\sum\limits_{j=1}^{c(G^d)}b_j\alpha_j({\rm
mod} 2). \ \ \ (2.4)$$

If there is another decomposition

$$\widetilde{C}\equiv\sum\limits_{M\in V(G^d)}
\sum\limits_{i=1}^{l'(M)}a_{i}^{'M}C_{i}^M+\sum\limits_{j=1}^{c'(G^d)}b'_{j}\alpha_{j}({\rm
mod} 2),$$

\no not loss of generality, assume $l'(M)\leq l(M)$ and $c'(M)\leq
c(M)$, then we know that

$$\sum\limits_{M\in V(G^d)}
\sum\limits_{i=1}^{l(M)}(a_i^M-a_i^{'M})C_i^M+
\sum\limits_{j=1}^{c(G^d)}(b_j-b'_j)\alpha_{j'}=0,$$

\no where ${a'}_i^M=0$ if $i>l'(M)$, $b'_j=0$ if $j'>c'(M)$. Since
$C_i^M, 1\leq i\leq b(M)$ and $\alpha_j, 1\leq j\leq\beta(G^d)$ are
bases of the fundamental group $\pi(M)$ and $\pi(G^d)$ respectively,
we must have

$$a_i^M=a_i^{'M}, 1\leq i\leq l(M) \ {\rm and} \ b_j=b'_j, 1\leq j\leq c(G^d).$$

\no Whence, the decomposition $(2.4)$ is unique.

For proving the conclusion ($2$), notice that
$\widetilde{M}(n_1,n_2,\cdots,n_m)$ is pathwise $d$-connected. Let
$P^d(x,y)$ be a $d$-path connecting points $x$ and $y$ in
$\widetilde{M}(n_1,n_2,\cdots,n_m)$. Define

$$\omega_*(C)=P^d(x,y)C(P^d)^{-1}(x,y)$$

\no for $\forall C\in\widetilde{M}(n_1,n_2,\cdots,n_m)$. Then it
can be checked immediately that

$$\omega_*:\pi^{d}(\widetilde{M}(n_1,n_2,\cdots,n_m),x)\rightarrow
\pi^{d}(\widetilde{M}(n_1,n_2,\cdots,n_m),y)$$

\no is an isomorphism.\ \  \ \ $\natural$

A $d$-connected finitely combinatorial manifold
$\widetilde{M}(n_1,n_2,\cdots,n_m)$ is said to be {\it simply
$d$-connected} if $\pi^d(\widetilde{M}(n_1,n_2,\cdots,n_m),x)$ is
trivial. As a consequence, we get the following result by Theorem
$2.7$.

\vskip 4mm

\no{\bf Corollary $2.3$} \ {\it A $d$-connected finitely
combinatorial manifold $\widetilde{M}(n_1,n_2,\cdots,n_m)$ is simply
$d$-connected if and only if}

($1$) \ {\it for $\forall M\in
V(G^d[\widetilde{M}(n_1,n_2,\cdots,n_m)])$, $M$ is simply
$d$-connected and}

($2$) \ {\it $G^d[\widetilde{M}(n_1,n_2,\cdots,n_m)]$ is a tree.}

\vskip 3mm

{\it Proof} \ According to the decomposition for
$\pi^d(\widetilde{M}(n_1,n_2,\cdots,n_m),x)$ in Theorem $2.7$, it
is trivial if and only if $\pi(M)$ and $\pi(G^d)$ both are trivial
for $\forall M\in V(G^d[\widetilde{M}(n_1,n_2,\cdots,n_m)])$, i.e
$M$ is simply $d$-connected and $G^d$ is a tree. \ \  \ $\natural$

For equivalent homotopically combinatorial manifolds, we also get a
criterion under a homotopically equivalent mapping in the next.

\vskip 4mm

\no{\bf Theorem $2.8$} \ {\it If
$f:\widetilde{M}(n_1,n_2,\cdots,n_m)\rightarrow
\widetilde{M}(k_1,k_2,\cdots,k_l)$ is a homotopic equivalence, then
for any integer $d, 1\leq d\leq n_1 $ and
$x\in\widetilde{M}(n_1,n_2,\cdots,n_m)$,}

$$\pi^d(\widetilde{M}(n_1,n_2,\cdots,n_m),x)\cong
\pi^d(\widetilde{M}(k_1,k_2,\cdots,k_l),f(x)).$$

\vskip 3mm

{\it Proof} \ Notice that $f$ can natural induce a homomorphism

$$f_{\pi}:\pi^d(\widetilde{M}(n_1,n_2,\cdots,n_m),x)\rightarrow
\pi^d(\widetilde{M}(k_1,k_2,\cdots,k_l),f(x))$$

\no defined by $f_{\pi}\left<g\right>=\left<f(g)\right>$ for
$\forall g\in\pi^d(\widetilde{M}(n_1,n_2,\cdots,n_m),x)$ since it
can be easily checked that $f_{\pi}(gh)=f_{\pi}(g)f_{\pi}(h)$ for
$\forall g,h\in\pi^d(\widetilde{M}(n_1,n_2,\cdots,n_m),x)$. We
only need to prove that $f_{\pi}$ is an isomorphism.

By definition, there is also a homotopic equivalence
$g:\widetilde{M}(k_1,k_2,\cdots,k_l)\rightarrow\widetilde{M}(n_1,n_2,\cdots,n_m)$
such that $gf\simeq
identity:\widetilde{M}(n_1,n_2,\cdots,n_m)\rightarrow\widetilde{M}(n_1,n_2,\cdots,n_m)$.
Thereby, $g_{\pi}f_{\pi}=(gf)_{\pi}=\mu(identity)_{\pi}:$

$$\pi^d(\widetilde{M}(n_1,n_2,\cdots,n_m),x)
\rightarrow\pi^s(\widetilde{M}(n_1,n_2,\cdots,n_m),x),$$

\no where $\mu$ is an isomorphism induced by a certain $d$-path
from $x$ to $gf(x)$ in $\widetilde{M}(n_1,n_2,$ $\cdots,n_m)$.
Therefore, $g_{\pi}f_{\pi}$ is an isomorphism. Whence, $f_{\pi}$
is a monomorphism and $g_{\pi}$ is an epimorphism.

Similarly, apply the same argument to the homotopy

$$fg\simeq
identity:\widetilde{M}(k_1,k_2,\cdots,k_l)\rightarrow\widetilde{M}(k_1,k_2,\cdots,k_l),$$

\no we get that
$f_{\pi}g_{\pi}=(fg)_{\pi}=\nu(identity)_{pi}:$

$$\pi^d(\widetilde{M}(k_1,k_2,\cdots,k_l),x)
\rightarrow\pi^s(\widetilde{M}(k_1,k_2,\cdots,k_l),x),$$

\no where $\nu$ is an isomorphism induced by a $d$-path from
$fg(x)$ to $x$ in $\widetilde{M}(k_1,k_2,$ $\cdots,k_l)$. So
$g_{\pi}$ is a monomorphism and $f_{\pi}$ is an epimorphism.
Combining these facts enables us to conclude that
$f_{\pi}:\pi^d(\widetilde{M}(n_1,n_2,\cdots,n_m),x)\rightarrow
\pi^d(\widetilde{M}(k_1,k_2,\cdots,k_l),f(x))$ is an isomorphism
.\ \ $\natural$

\vskip 4mm

\no{\bf Corollary $2.4$} \ {If
$f:\widetilde{M}(n_1,n_2,\cdots,n_m)\rightarrow
\widetilde{M}(k_1,k_2,\cdots,k_l)$ is a homeomorphism, then for any
integer $d, 1\leq d\leq n_1 $ and
$x\in\widetilde{M}(n_1,n_2,\cdots,n_m)$,}

$$\pi^d(\widetilde{M}(n_1,n_2,\cdots,n_m),x)\cong
\pi^d(\widetilde{M}(k_1,k_2,\cdots,k_l),f(x)).$$

\vskip 5mm

\no{\bf $2.4$ Euler-Poincare characteristic}

\vskip 3mm

\no It is well-known that the integer

$$\chi(\mathfrak{M})=\sum\limits_{i=0}^{\infty}(-1)^i\alpha_i$$

\no with $\alpha_i$ the number of $i$-dimensional cells in a
$CW$-complex $\mathfrak{M}$ is defined to be the Euler-Poincare
characteristic of this complex. In this subsection, we get the
Euler-Poincare characteristic for finitely combinatorial manifolds.
For this objective, define a clique sequence $\{Cl(i)\}_{i\geq 1}$
in the graph $G[\widetilde{M}]$ by the following programming.

\vskip 2mm

\no{STEP $1.$} Let $Cl(G[\widetilde{M}])=l_0$. Construct

\begin{eqnarray*}
Cl(l_0)&=&\{K_1^{l_0},K_2^{l_0},\cdots,K_{p}^{i_0}| K_i^{l_0}\succ
G[\widetilde{M}] \ {\rm and} \ K_i^{l_0}\cap K_j^{l_0} =\emptyset, \\
& \ & {\rm or \ a \ vertex\in V(G[\widetilde{M}]) \ for}\
i\not=j,1\leq i,j\leq p\}.
\end{eqnarray*}

\no{STEP $2.$} \ Let $G_1=\bigcup\limits_{K^l\in Cl(l)}K^l$ and
$Cl(G[\widetilde{M}]\setminus G_1)=l_1$. Construct

\begin{eqnarray*}
Cl(l_1)&=&\{K_1^{l_1},K_2^{l_1},\cdots,K_{q}^{i_1}| K_i^{l_1}\succ
G[\widetilde{M}] \ {\rm and} \ K_i^{l_1}\cap K_j^{l_1} =\emptyset\\
& \ & {\rm or \ a \ vertex\in V(G[\widetilde{M}]) \ for}\
i\not=j,1\leq i,j\leq q\}.
\end{eqnarray*}

\no{STEP $3.$} Assume we have constructed $Cl(l_{k-1})$ for an
integer $k\geq 1$. Let $G_k=\bigcup\limits_{K^{l_{k-1}}\in
Cl(l)}K^{l_{k-1}}$ and
$Cl(G[\widetilde{M}]\setminus(G_1\cup\cdots\cup G_k))=l_k$. We
construct

\begin{eqnarray*}
Cl(l_k)&=&\{K_1^{l_k},K_2^{l_k},\cdots,K_{r}^{l_k}| K_i^{l_k}\succ
G[\widetilde{M}] \ {\rm and} \ K_i^{l_k}\cap K_j^{l_k} =\emptyset,\\
& \ & {\rm or \ a \ vertex\in V(G[\widetilde{M}]) \ for}\
i\not=j,1\leq i,j\leq r\}.
\end{eqnarray*}

\no{STEP $4.$} Continue STEP $3$ until we find an integer $t$ such
that there are no edges in
$G[\widetilde{M}]\setminus\bigcup\limits_{i=1}^tG_i$.\vskip 2mm

By this clique sequence $\{Cl(i)\}_{i\geq 1}$, we can calculate the
Eucler-Poincare characteristic of finitely combinatorial manifolds.

\vskip 4mm

\no{\bf Theorem $2.9$} \ {\it Let $\widetilde{M}$ be a finitely
combinatorial manifold. Then}

$$\chi(\widetilde{M})=\sum\limits_{K^k\in Cl(k), k\geq 2}
\sum\limits_{M_{i_j}\in V(K^k),1\leq j\leq s\leq
k}(-1)^{s+1}\chi(M_{i_1}\bigcup\cdots\bigcup M_{i_s})$$

\vskip 3mm

{\it Proof} \ Denoted  the numbers of all these $i$-dimensional
cells in a combinatorial manifold $\widetilde{M}$ or in a manifold
$M$ by $\widetilde{\alpha}_i$ and $\alpha_i(M)$. If
$G[\widetilde{M}]$ is nothing but a complete graph $K^k$ with
$V(G[\widetilde{M}])=\{M_1,M_2,\cdots,M_k\}$, $k\geq 2$, by applying
the inclusion-exclusion principe and the definition of
Euler-Poincare characteristic we get that

\begin{eqnarray*}\chi(\widetilde{M})&=&
\sum\limits_{i=0}^{\infty}(-1)^i\widetilde{\alpha}_i\\
&=& \sum\limits_{i=0}^{\infty}(-1)^i\sum\limits_{M_{i_j}\in
V(K^k),1\leq j\leq s\leq
k}(-1)^{s+1}\alpha_i(M_{i_1}\bigcup\cdots\bigcup
M_{i_s})\\
&=& \sum\limits_{M_{i_j}\in V(K^k),1\leq j\leq s\leq
k}(-1)^{s+1}\sum\limits_{i=0}^{\infty}(-1)^i\alpha_i(M_{i_1}\bigcup\cdots\bigcup
M_{i_s})\\
&=& \sum\limits_{M_{i_j}\in V(K^k),1\leq j\leq s\leq
k}(-1)^{s+1}\chi(M_{i_1}\bigcup\cdots\bigcup M_{i_s})\end{eqnarray*}

\no for instance,
$\chi(\widetilde{M})=\chi(M_1)+\chi(M_2)-\chi(M_1\cap M_2)$ if
$G[\widetilde{M}]=K^2$ and $V(G[\widetilde{M}])=\{M_1,M_2\}$. By the
definition of clique sequence of $G[\widetilde{M}]$, we finally
obtain that

$$\chi(\widetilde{M})=\sum\limits_{K^k\in Cl(k),k\geq 2}
\sum\limits_{M_{i_j}\in V(K^k),1\leq j\leq s\leq
k}(-1)^{i+1}\chi(M_{i_1}\bigcup\cdots\bigcup M_{i_s}).  \ \ \
\natural$$

If $G[\widetilde{M}]$ is just one of some special graphs, we can get
interesting consequences by Theorem $2.9$.

\vskip 4mm

\no{\bf Corollary $2.5$} \ {\it Let $\widetilde{M}$ be a finitely
combinatorial manifold. If $G[\widetilde{M}]$ is $K^3$-free, then}

$$\chi(\widetilde{M})=\sum\limits_{M\in V(G[\widetilde{M}])}\chi^2(M)
-\sum\limits_{(M_1,M_2)\in E(G[\widetilde{M}])}\chi(M_1\bigcap
M_2).$$

{\it Particularly, if ${\rm dim}(M_1\bigcap M_2)$ is a constant for
any $(M_1,M_2)\in E(G[\widetilde{M}])$, then}

$$\chi(\widetilde{M})=\sum\limits_{M\in V(G[\widetilde{M}])}\chi^2(M)
-\chi(M_1\bigcap M_2)|E(G[\widetilde{M}])|.$$

\vskip 3mm

{\it Proof} \ Notice that $G[\widetilde{M}]$ is $K^3$-free, we get
that

\begin{eqnarray*}
\chi(\widetilde{M})&=& \sum\limits_{(M_1,M_2)\in
E(G[\widetilde{M}])}(\chi(M_1)+\chi(M_2)-\chi(M_1\bigcap M_2))\\
&=& \sum\limits_{(M_1,M_2)\in
E(G[\widetilde{M}])}(\chi(M_1)+\chi(M_2))+\sum\limits_{(M_1,M_2)\in
E(G[\widetilde{M}])}\chi(M_1\bigcap M_2))\\
&=& \sum\limits_{M\in V(G[\widetilde{M}])}\chi^2(M)
-\sum\limits_{(M_1,M_2)\in E(G[\widetilde{M}])}\chi(M_1\bigcap M_2).
\ \ \ \natural
\end{eqnarray*}

Since the Euler-Poincare characteristic of a manifold $M$ is $0$ if
${\rm dim}M\equiv 1(mod 2)$, we get the following consequence.

\vskip 4mm

\no{\bf Corollary $2.6$} \ {\it Let $\widetilde{M}$ be a finitely
combinatorial manifold with odd dimension number for any
intersection of $k$ manifolds with $k\geq 2$. Then}

$$\chi(\widetilde{M})=\sum\limits_{M\in V(G[\widetilde{M}])}\chi(M).$$

\vskip 6mm

\no{\bf \S $3.$ Differential structures on combinatorial
manifolds}

\vskip 4mm

\no We introduce differential structures on finitely combinatorial
manifolds and characterize them in this section.

\vskip 4mm

\no{\bf $3.1$ Tangent vector fields}

\vskip 3mm

\no{\bf Definition $3.1$} \ {\it For a given integer sequence $1\leq
n_1<n_2<\cdots<n_m$, a combinatorially $C^h$ differential manifold
$(\widetilde{M}(n_1,n_2,\cdots,n_m); \widetilde{{\mathcal A}})$ is a
finitely combinatorial manifold $\widetilde{M}(n_1,n_2,\cdots,n_m)$,
$\widetilde{M}(n_1,n_2,\cdots,n_m)=\bigcup\limits_{i\in I}U_i$,
endowed with a atlas $\widetilde{{\mathcal
A}}=\{(U_{\alpha};\varphi_{\alpha})| \alpha\in I\}$ on
$\widetilde{M}(n_1,n_2,\cdots,n_m)$ for an integer $h, h\geq 1$ with
conditions following hold.

$(1)$ \ $\{U_{\alpha}; \alpha\in I\}$ is an open covering of
$\widetilde{M}(n_1,n_2,\cdots,n_m)$;

$(2)$ \ For $\forall \alpha,\beta\in I$, local charts
$(U_{\alpha};\varphi_{\alpha})$ and $(U_{\beta};\varphi_{\beta})$
are {\it equivalent}, i.e., $U_{\alpha}\bigcap U_{\beta}=\emptyset$
or $U_{\alpha}\bigcap U_{\beta}\not=\emptyset$ but the {\it overlap
maps}

$$
\varphi_{\alpha}\varphi_{\beta}^{-1}:
\varphi_{\beta}(U_{\alpha}\bigcap U_{\beta})\rightarrow
\varphi_{\beta}(U_{\beta}) \ \ {\rm and} \ \
\varphi_{\beta}\varphi_{\alpha}^{-1}:
\varphi_{\beta}(U_{\alpha}\bigcap U_{\beta})\rightarrow
\varphi_{\alpha}(U_{\alpha})
$$

\no are $C^h$ mappings;

$(3)$ \ $\widetilde{{\mathcal A}}$ is maximal, i.e., if
$(U;\varphi)$ is a local chart of
$\widetilde{M}(n_1,n_2,\cdots,n_m)$ equivalent with one of local
charts in $\widetilde{{\mathcal A}}$, then
$(U;\varphi)\in\widetilde{{\mathcal A}}$.

Denote by $(\widetilde{M}(n_1,n_2,\cdots,n_m);\widetilde{{\mathcal
A}})$ a combinatorially differential manifold. A finitely
combinatorial manifold $\widetilde{M}(n_1,n_2,\cdots,n_m)$ is said
to be {\it smooth} if it is endowed with a $C^{\infty}$ differential
structure.}

\vskip 3mm

Let $\widetilde{{\mathcal A}}$ be an atlas on
$\widetilde{M}(n_1,n_2,\cdots,n_m)$. Choose a local chart
$(U;\varpi)$ in $\widetilde{{\mathcal A}}$. For $\forall
p\in(U;\varphi)$, if $\varpi_p:
U_p\rightarrow\bigcup\limits_{i=1}^{s(p)}B^{n_i(p)}$ and
$\widehat{s}(p)= {\rm dim}(\bigcap\limits_{i=1}^{s(p)}B^{n_i(p)})$,
the following $s(p)\times n_{s(p)}$ matrix $[\varpi(p)]$

\[
[\varpi(p)]=\left[
\begin{array}{cccccccc}
\frac{x^{11}}{\widehat{s}(p)} & \cdots & \frac{x^{1\widehat{s}(p)}}{\widehat{s}(p)}
& x^{1(\widehat{s}(p)+1)} & \cdots & x^{1n_1} & \cdots & 0 \\
\frac{x^{21}}{\widehat{s}(p)} & \cdots &
\frac{x^{2\widehat{s}(p)}}{\widehat{s}(p)}
& x^{2(\widehat{s}(p)+1)} & \cdots & x^{2n_2} & \cdots & 0  \\
\cdots & \cdots & \cdots  & \cdots & \cdots & \cdots  \\
\frac{x^{s(p)1}}{\widehat{s}(p)} & \cdots &
\frac{x^{s(p)\widehat{s}(p)}}{\widehat{s}(p)} &
x^{s(p)(\widehat{s}(p)+1)} & \cdots & \cdots & x^{s(p)n_{s(p)}-1} &
x^{s(p)n_{s(p)}}
\end{array}
\right]
\]

\no with $x^{is}=x^{js}$ for $1\leq i,j\leq s(p), 1\leq
s\leq\widehat{s}(p)$ is called the {\it coordinate matrix of $p$}.
For emphasize $\varpi$ is a matrix, we often denote local charts in
a combinatorially differential manifold by $(U;[\varpi])$. Using the
coordinate matrix system of a combinatorially differential manifold
$(\widetilde{M}(n_1,n_2,\cdots,n_m); \widetilde{{\mathcal A}})$, we
introduce the conception of $C^h$ mappings and functions in the
next.

\vskip 4mm

\no{\bf Definition $3.2$} \ {\it Let
$\widetilde{M}_1(n_1,n_2,\cdots,n_m)$,
$\widetilde{M}_2(k_1,k_2,\cdots,k_l)$ be smoothly combinatorial
manifolds and}

$$f:\widetilde{M}_1(n_1,n_2,\cdots,n_m)\rightarrow\widetilde{M}_2(k_1,k_2,\cdots,k_l)$$

\no{\it be a mapping, $p\in\widetilde{M}_1(n_1,n_2,\cdots,n_m)$. If
there are local charts $(U_p;[\varpi_p])$ of $p$ on
$\widetilde{M}_1(n_1,n_2,\cdots,n_m)$ and
$(V_{f(p)};[\omega_{f(p)}])$ of $f(p)$ with $f(U_p)\subset V_{f(p)}$
such that the composition mapping}

$$\widetilde{f}=[\omega_{f(p)}]\circ f\circ[\varpi_p]^{-1}:
[\varpi_p](U_p)\rightarrow[\omega_{f(p)}](V_{f(p)})$$

\no{is a $C^h$ mapping, then $f$ is called a $C^h$ mapping at the
point $p$. If $f$ is $C^h$ at any point $p$ of
$\widetilde{M}_1(n_1,n_2,\cdots,n_m)$, then $f$ is called a $C^h$
mapping. Particularly, if $\widetilde{M}_2(k_1,k_2,\cdots,k_l)={\bf
R}$, $f$ ia called a $C^h$ function on
$\widetilde{M}_1(n_1,n_2,\cdots,n_m)$. In the extreme $h=\infty$,
these terminologies are called smooth mappings and functions,
respectively. Denote by $\mathscr{X}_p$ all these $C^{\infty}$
functions at a point $p\in\widetilde{M}(n_1,n_2,\cdots,n_m)$.}

\vskip 3mm

For the existence of combinatorially differential manifolds, we
know the following result.

\vskip 4mm

\no{\bf Theorem $3.1$} \ {\it Let
$\widetilde{M}(n_1,n_2,\cdots,n_m)$ be a finitely combinatorial
manifold and $d,1\leq d\leq n_1$ an integer. If $\forall M\in
V(G^d[\widetilde{M}(n_1,n_2,\cdots,n_m)])$ is $C^h$ differential and
$\forall (M_1,M_2)\in E(G^d[\widetilde{M}(n_1,n_2,\cdots,n_m)])$
there exist atlas}

$${\mathcal A}_1=\{(V_x;\varphi_x)| \forall x\in M_1\} \ \
{\mathcal A}_2=\{(W_y;\psi_y)| \forall y\in M_2\}$$

\no{\it such that $\varphi_x|_{V_x\bigcap W_y}=\psi_y|_{V_x\bigcap
W_y}$ for $\forall x\in M_1, y\in M_2$, then there is a
differential structures}

$$\widetilde{{\mathcal A}}=\{(U_p;[\varpi_p])| \forall
p\in\widetilde{M}(n_1,n_2,\cdots,n_m)\}$$

\no{\it such that
$(\widetilde{M}(n_1,n_2,\cdots,n_m);\widetilde{{\mathcal A}})$ is a
combinatorially $C^h$ differential manifold.}

\vskip 3mm

{\it Proof} \ By definition, We only need to show that we can
always choose a neighborhood $U_p$ and a homoeomorphism
$[\varpi_p]$ for each $p\in\widetilde{M}(n_1,n_2,\cdots,n_m)$
satisfying these conditions $(1)-(3)$ in definition $3.1$.

By assumption, each manifold $\forall M\in
V(G^d[\widetilde{M}(n_1,n_2,\cdots,n_m)])$ is $C^h$ differential,
accordingly there is an index set $I_M$ such that $\{U_{\alpha};
\alpha\in I_M\}$ is an open covering of $M$, local charts
$(U_{\alpha};\varphi_{\alpha})$ and $(U_{\beta};\varphi_{\beta})$ of
$M$ are equivalent and ${\mathcal A}=\{(U;\varphi)\}$ is maximal.
Since for $\forall p\in\widetilde{M}(n_1,n_2,\cdots,n_m)$, there is
a local chart $(U_p;[\varpi_p])$ of $p$ such that $[\varpi_p]:
U_p\rightarrow\bigcup\limits_{i=1}^{s(p)}B^{n_i(p)}$, i.e., $p$ is
an intersection point of manifolds $M^{n_i(p)},1\leq i\leq s(p)$. By
assumption each manifold $M^{n_i(p)}$ is $C^h$ differential, there
exists a local chart $(V_p^i; \varphi_p^i)$ while the point $p\in
M^{n_i(p)}$ such that $\varphi_p^i\rightarrow B^{n_i(p)}$. Now we
define

$$U_p=\bigcup\limits_{i=1}^{s(p)}V_p^i.$$

\no Then applying the Gluing Lemma again, we know that there is a
homoeomorphism $[\varpi_p]$ on $U_p$ such that

$$[\varpi_p]|_{M^{n_i(p)}}=\varphi_p^i$$

\no for any integer $i, \leq i\leq s(p)$. Thereafter,

$$\widetilde{{\mathcal A}}=\{(U_p;[\varpi_p])| \forall
p\in\widetilde{M}(n_1,n_2,\cdots,n_m)\}$$

\no is a $C^h$ differential structure on
$\widetilde{M}(n_1,n_2,\cdots,n_m)$ satisfying conditions $(1)-(3)$.
Thereby $(\widetilde{M}(n_1,n_2,\cdots,n_m);\widetilde{{\mathcal
A}})$ is a combinatorially $C^h$ differential manifold. \ \ \
$\natural$

\vskip 4mm

\no{\bf Definition $3.3$} \ {\it Let
$(\widetilde{M}(n_1,n_2,\cdots,n_m),\widetilde{{\mathcal A}})$ be a
smoothly combinatorial manifold and
$p\in\widetilde{M}(n_1,n_2,\cdots,n_m)$. A tangent vector $v$ at $p$
is a mapping $v: \mathscr{X}_p\rightarrow {\bf R}$ with conditions
following hold.}\vskip 2mm

$(1)$ \ $\forall g,h\in\mathscr{X}_p, \forall\lambda\in {\bf R}, \
v(h+\lambda h)=v(g)+\lambda v(h);$

$(2)$ \ $\forall g,h\in\mathscr{X}_p, v(gh)=v(g)h(p)+g(p)v(h).$

\vskip 3mm

Denoted all tangent vectors at
$p\in\widetilde{M}(n_1,n_2,\cdots,n_m)$ by
$T_p\widetilde{M}(n_1,n_2,\cdots,n_m)$ and define addition¡°+¡±and
scalar multiplication¡°$\cdot$¡±for $\forall u,v\in
T_p\widetilde{M}(n_1,n_2,\cdots,n_m),$ $\lambda\in{\bf R}$ and
$f\in\mathscr{X}_p$ by

$$(u+v)(f)=u(f)+v(f), \ \ (\lambda u)(f)=\lambda\cdot u(f).$$

\no Then it can be shown immediately that
$T_p\widetilde{M}(n_1,n_2,\cdots,n_m)$ is a vector space under
these two operations¡°+¡±and¡°$\cdot$¡±.

\vskip 4mm

\no{\bf Theorem $3.2$} \ {\it For any point
$p\in\widetilde{M}(n_1,n_2,\cdots,n_m)$ with a local chart $(U_p;
[\varphi_p])$, the dimension of
$T_p\widetilde{M}(n_1,n_2,\cdots,n_m)$ is}

$${\rm dim}T_p\widetilde{M}(n_1,n_2,\cdots,n_m) =\widehat{s}(p)+ \sum\limits_{i=1}^{s(p)}(n_i-\widehat{s}(p))$$

\no{\it with a basis matrix }\vskip 3mm

$$[\frac{\partial}{\partial\overline{x}}]_{s(p)\times n_{s(p)}}=\hskip 115mm$$

\[
\left[
\begin{array}{cccccccc}
\frac{1}{\widehat{s}(p)}\frac{\partial}{\partial x^{11}} & \cdots &
\frac{1}{\widehat{s}(p)}\frac{\partial}{\partial
x^{1\widehat{s}(p)}}
& \frac{\partial}{\partial x^{1(\widehat{s}(p)+1)}} & \cdots & \frac{\partial}{\partial x^{1n_1}} & \cdots & 0 \\
\frac{1}{\widehat{s}(p)}\frac{\partial}{\partial x^{21}} & \cdots &
\frac{1}{\widehat{s}(p)}\frac{\partial}{\partial
x^{2\widehat{s}(p)}}
& \frac{\partial}{\partial x^{2(\widehat{s}(p)+1)}} & \cdots & \frac{\partial}{\partial x^{2n_2}} & \cdots & 0 \\
\cdots & \cdots & \cdots  & \cdots & \cdots & \cdots  \\
\frac{1}{\widehat{s}(p)}\frac{\partial}{\partial x^{s(p)1}} & \cdots
& \frac{1}{\widehat{s}(p)}\frac{\partial}{\partial
x^{s(p)\widehat{s}(p)}} & \frac{\partial}{\partial
x^{s(p)(\widehat{s}(p)+1)}} & \cdots & \cdots &
\frac{\partial}{\partial x^{s(p)(n_{s(p)}-1)}} &
\frac{\partial}{\partial x^{s(p)n_{s(p)}}}
\end{array} \right]
\]

\vskip 3mm

\no{\it where $x^{il}=x^{jl}$ for $1\leq i,j\leq s(p), 1\leq
l\leq\widehat{s}(p)$, namely there is a smoothly functional matrix
$[v_{ij}]_{s(p)\times n_{s(p)}}$ such that for any tangent vector
$\overline{v}$ at a point $p$ of
$\widetilde{M}(n_1,n_2,\cdots,n_m)$,}

$$\overline{v}=[v_{ij}]_{s(p)\times n_{s(p)}}\odot[\frac{\partial}
{\partial\overline{x}}]_{s(p)\times n_{s(p)}},$$

\no{\it where $[a_{ij}]_{k\times l}\odot[b_{ts}]_{k\times
l}=\sum\limits_{i=1}^k\sum\limits_{j=1}^la_{ij}b_{ij}$.}

\vskip 3mm

{\it Proof} \ For $\forall f\in\mathscr{X}_p$, let
$\widetilde{f}=f\cdot[\varphi_p]^{-1}\in
\mathscr{X}_{[\varphi_p](p)}$. We only need to prove that $f$ can be
spanned by elements in

$$\{\frac{\partial}{\partial x^{hj}}|_p|1\leq j\leq\widehat{s}(p)\}\bigcup
(\bigcup\limits_{i=1}^{s(p)}\bigcup\limits_{j=\widehat{s}(p)+1}^{n_i}
\{\frac{\partial}{\partial x^{ij}}|_p \ | \ 1\leq j\leq s\}), \ \
(3.1)$$

\no for a given integer $h,1\leq h\leq s(p)$, namely $(3.1)$ is a
basis of $T_p\widetilde{M}(n_1,n_2,\cdots,n_m)$. In fact, for
$\forall\overline{x}\in[\varphi_p](U_p)$, since $\widetilde{f}$ is
smooth, we know that

\begin{eqnarray*}
\widetilde{f}(\overline{x})-\widetilde{f}(\overline{x}_0)&=&
\int\limits_{0}^{1}\frac{d}{dt}\widetilde{f}(\overline{x}_0+t(\overline{x}-\overline{x}_0))dt\\
&=&
\sum\limits_{i=1}^{s(p)}\sum\limits_{j=1}^{n_i}\eta_{\widehat{s}(p)}^j(x^{ij}-x_0^{ij})
\int\limits_{0}^1\frac{\partial\widetilde{f}}{\partial
x^{ij}}(\overline{x}_0+t(\overline{x}-\overline{x}_0))dt
\end{eqnarray*}

\no in a spherical neighborhood of the point $p$ in
$[\varphi_p](U_p)\subset{\bf
R}^{\widehat{s}(p)-s(p)\widehat{s}(p)+n_1+n_2+\cdots+n_{s(p)}}$ with
$[\varphi_p](p)=\overline{x}_0$, where

\[
\eta_{\widehat{s}(p)}^j=\left\{\begin{array}{cc}
\frac{1}{\widehat{s}(p)},& {\rm if}\quad 1\leq j\leq\widehat{s}(p) ,\\
1,& {\rm otherwise}.
\end{array}
\right.
\]

\no Define

$$\widetilde{g}_{ij}(\overline{x})=\int\limits_0^1\frac{\partial\widetilde{f}}{\partial
x^{ij}}(\overline{x}_0+t(\overline{x}-\overline{x}_0))dt$$

\no and $g_{ij}=\widetilde{g}_{ij}\cdot[\varphi_p]$. Then we find
that

\begin{eqnarray*}
g_{ij}(p)=\widetilde{g}_{ij}(\overline{x}_0)&=&
\frac{\partial\widetilde{f}}{\partial x^{ij}}(\overline{x}_0)\\
&=& \frac{\partial(f\cdot[\varphi_p]^{-1})}{\partial
x^{ij}}([\varphi_p](p))=\frac{\partial f}{\partial x^{ij}}(p).
\end{eqnarray*}

\no Therefore, for $\forall q\in U_p$, there are $g_{ij}, 1\leq
i\leq s(p), 1\leq j\leq n_i$ such that

$$f(q)=f(p)+
\sum\limits_{i=1}^{s(p)}\sum\limits_{j=1}^{n_i}\eta_{\widehat{s}(p)}^j(x^{ij}-x_0^{ij})g_{ij}(p).$$

Now let $\overline{v}\in T_p\widetilde{M}(n_1,n_2,\cdots,n_m)$.
Application of the condition $(2)$ in Definition $3.1$ shows that

$$v(f(p))=0, \ \ {\rm and} \ \ v(\eta_{\widehat{s}(p)}^jx_0^{ij})=0.$$

\no Accordingly, we obtain that

\begin{eqnarray*}
\overline{v}(f)&=& \overline{v}(f(p)+\sum\limits_{i=1}^{s(p)}
\sum\limits_{j=1}^{n_i}\eta_{\widehat{s}(p)}^j(x^{ij}-x_0^{ij})g_{ij}(p))\\
&=& \overline{v}(f(p))+\sum\limits_{i=1}^{s(p)}
\sum\limits_{j=1}^{n_i}\overline{v}(\eta_{\widehat{s}(p)}^j(x^{ij}-x_0^{ij})g_{ij}(p)))\\
&=& \sum\limits_{i=1}^{s(p)}
\sum\limits_{j=1}^{n_i}(\eta_{\widehat{s}(p)}^jg_{ij}(p)\overline{v}(x^{ij}-x_0^{ij})+(x^{ij}(p)-x_0^{ij})
\overline{v}(\eta_{\widehat{s}(p)}^jg_{ij}(p)))\\
&=& \sum\limits_{i=1}^{s(p)}
\sum\limits_{j=1}^{n_i}\eta_{\widehat{s}(p)}^j\frac{\partial
f}{\partial
x^{ij}}(p)\overline{v}(x^{ij})\\
&=& \sum\limits_{i=1}^{s(p)}
\sum\limits_{j=1}^{n_i}\overline{v}(x^{ij})\eta_{\widehat{s}(p)}^j\frac{\partial}{\partial
x^{ij}}|_p(f)=[v_{ij}]_{s(p)\times n_{s(p)}}\odot[\frac{\partial}
{\partial\overline{x}}]_{s(p)\times n_{s(p)}}|_p(f).
\end{eqnarray*}

\no Therefore, we get that

$$\overline{v}=[v_{ij}]_{s(p)\times n_{s(p)}}\odot[\frac{\partial}
{\partial\overline{x}}]_{s(p)\times n_{s(p)}}. \ \ (3.2)$$

The formula $(3.2)$ shows that any tangent vector $\overline{v}$
in $T_p\widetilde{M}(n_1,n_2,\cdots,n_m)$ can be spanned by
elements in $(3.1)$.

Notice that all elements in $(3.1)$ are also linearly independent.
Otherwise, if there are numbers $a^{ij}, 1\leq i\leq s(p), 1\leq
j\leq n_i$ such that

$$(\sum\limits_{j=1}^{\widehat{s}(p)}a^{hj}\frac{\partial}{\partial x^{hj}}+
\sum\limits_{i=1}^{s(p)}\sum\limits_{j=\widehat{s}(p)+1}^{n_i}a^{ij}
\frac{\partial}{\partial x^{ij}})|_p=0,$$

\no then we get that

$$a^{ij}=(\sum\limits_{j=1}^{\widehat{s}(p)}a^{hj}\frac{\partial}{\partial x^{hj}}+
\sum\limits_{i=1}^{s(p)}\sum\limits_{j=\widehat{s}(p)+1}^{n_i}a^{ij}
\frac{\partial}{\partial x^{ij}})(x^{ij})=0$$

\no for $1\leq i\leq s(p), 1\leq j\leq n_i$. Therefore, $(3.1)$ is a
basis of the tangent vector space
$T_p\widetilde{M}(n_1,n_2,\cdots,n_m)$ at the point
$p\in(\widetilde{M}(n_1,n_2,\cdots,n_m);\widetilde{{\mathcal A}})$.
\ \ $\natural$\vskip 2mm

By Theorem $3.2$, if $s(p)=1$ for any point
$p\in(\widetilde{M}(n_1,n_2,\cdots,n_m);\widetilde{{\mathcal A}})$,
then ${\rm dim}T_p\widetilde{M}(n_1,n_2,\cdots,n_m)=n_1$. This can
only happens while $\widetilde{M}(n_1,n_2,\cdots,n_m)$ is combined
by one manifold.  As a consequence, we get a well-known result in
classical differential geometry again.

\vskip 4mm

\no{\bf Corollary $3.1$}([2]) \ {\it Let $(M^n;{\mathcal A})$ be a
smooth manifold and $p\in M^n$. Then}

$${\rm dim}T_pM^n=n$$

\no{\it with a basis}

$$\{\frac{\partial}{\partial x^i}|_p \ | \ 1\leq i\leq n\}.$$

\vskip 3mm

\no{\bf Definition $3.4$} \ {\it For $\forall
p\in(\widetilde{M}(n_1,n_2,\cdots,n_m);\widetilde{{\mathcal A}})$,
the dual space $T_p^*\widetilde{M}(n_1,n_2,\cdots,n_m)$ is called a
co-tangent vector space at $p$.}

\vskip 4mm

\no{\bf Definition $3.5$} \ {\it For $f\in\mathscr{X}_p, d\in
T_p^*\widetilde{M}(n_1,n_2,\cdots,n_m)$ and $\overline{v}\in
T_p\widetilde{M}(n_1,n_2,\cdots,n_m)$, the action of $d$ on $f$,
called a differential operator $d: \mathscr{X}_p\rightarrow {\bf
R}$, is defined by}

$$df \ = \ \overline{v}(f).$$

\vskip 3mm

Then we immediately obtain the result following.

\vskip 4mm

\no{\bf Theorem $3.3$} \ {\it For $\forall
p\in(\widetilde{M}(n_1,n_2,\cdots,n_m);\widetilde{{\mathcal A}})$
with a local chart $(U_p; [\varphi_p])$, the dimension of
$T_p^*\widetilde{M}(n_1,n_2,\cdots,n_m)$ is}

$${\rm dim}T_p^*\widetilde{M}(n_1,n_2,\cdots,n_m) = \widehat{s}(p)+ \sum\limits_{i=1}^{s(p)}(n_i-\widehat{s}(p))$$

\no{\it with a basis matrix}

$$[d\overline{x}]_{s(p)\times n_{s(p)}}=\hskip 110mm$$

\[
\left[
\begin{array}{cccccccc}
\frac{dx^{11}}{\widehat{s}(p)} & \cdots &
\frac{dx^{1\widehat{s}(p)}}{\widehat{s}(p)}
& dx^{1(\widehat{s}(p)+1)} & \cdots & dx^{1n_1} & \cdots & 0 \\
\frac{dx^{21}}{\widehat{s}(p)} & \cdots &
\frac{dx^{2\widehat{s}(p)}}{\widehat{s}(p)}
& dx^{2(\widehat{s}(p)+1)} & \cdots & dx^{2n_2} & \cdots & 0  \\
\cdots & \cdots & \cdots  & \cdots & \cdots & \cdots  \\
\frac{dx^{s(p)1}}{\widehat{s}(p)} & \cdots &
\frac{dx^{s(p)\widehat{s}(p)}}{\widehat{s}(p)} &
dx^{s(p)(\widehat{s}(p)+1)} & \cdots & \cdots & dx^{s(p)n_{s(p)}-1}
& dx^{s(p)n_{s(p)}}
\end{array}
\right]
\]

\vskip 2mm

\no{\it where $x^{il}=x^{jl}$ for $1\leq i,j\leq s(p), 1\leq
l\leq\widehat{s}(p)$, namely for any co-tangent vector $d$ at a
point $p$ of $\widetilde{M}(n_1,n_2,\cdots,n_m)$, there is a
smoothly functional matrix $[u_{ij}]_{s(p)\times s(p)}$ such that,}

$$d=[u_{ij}]_{s(p)\times n_{s(p)}}\odot[d\overline{x}]_{s(p)\times n_{s(p)}}.$$

\vskip 5mm

\no{\bf $3.2$ Tensor fields}

\vskip 3mm

\no{\bf Definition $3.6$} \ {\it Let
$\widetilde{M}(n_1,n_2,\cdots,n_m)$ be a smoothly combinatorial
manifold and $p\in\widetilde{M}(n_1,n_2,\cdots,n_m)$. A tensor of
type $(r,s)$ at the point $p$ on
$\widetilde{M}(n_1,n_2,\cdots,n_m)$ is an $(r+s)$-multilinear
function $\tau$,}

$$\tau: \underbrace{T_p^*\widetilde{M}\times\cdots\times T_p^*\widetilde{M}}\limits_{r}
\times\underbrace{T_p\widetilde{M}\times\cdots\times
T_p\widetilde{M}}\limits_{s}\rightarrow{\bf R},$$

\no{\it where
$T_p\widetilde{M}=T_p\widetilde{M}(n_1,n_2,\cdots,n_m)$ and
$T_p^*\widetilde{M}=T_p^*\widetilde{M}(n_1,n_2,\cdots,n_m)$.}

\vskip 3mm

Denoted by $T_s^r(p,\widetilde{M})$ all tensors of type $(r,s)$ at
a point $p$ of $\widetilde{M}(n_1,n_2,\cdots,n_m)$. Then we know
its structure by Theorems $3.2$ and $3.3$.

\vskip 4mm

\no{\bf Theorem $3.4$} \ {\it  Let
$\widetilde{M}(n_1,n_2,\cdots,n_m)$ be a smoothly combinatorial
manifold and $p\in\widetilde{M}(n_1,n_2,\cdots,n_m)$. Then}

$$T_s^r(p,\widetilde{M})=\underbrace{T_p\widetilde{M}\otimes\cdots\otimes
T_p\widetilde{M}}\limits_{r}\otimes\underbrace{T_p^*\widetilde{M}\otimes\cdots\otimes
T_p^*\widetilde{M}}\limits_{s},$$

\no{\it where
$T_p\widetilde{M}=T_p\widetilde{M}(n_1,n_2,\cdots,n_m)$ and
$T_p^*\widetilde{M}=T_p^*\widetilde{M}(n_1,n_2,\cdots,n_m)$,
particularly, }

$${\rm dim}T_s^r(p,\widetilde{M})=(\widehat{s}(p)+ \sum\limits_{i=1}^{s(p)}(n_i-\widehat{s}(p)))^{r+s}.$$

\vskip 3mm

{\it Proof} \ By definition and multilinear algebra, any tensor
$t$ of type $(r,s)$ at the point $p$ can be uniquely written as

$$t = \sum t_{j_1\cdots j_s}^{i_1\cdots i_r}
\frac{\partial}{\partial
x^{i_1j_1}}|_p\otimes\cdots\otimes\frac{\partial}{\partial
x^{i_rj_r}}|_p \otimes dx^{k_1l_1}\otimes\cdots\otimes dx^{k_sl_s}$$

\no for components $ t_{j_1\cdots j_s}^{i_1\cdots i_r}\in{\bf R}$
according to Theorems $3.2$ and $3.3$, where $1\leq i_h, k_h\leq
s(p)$ and $1\leq j_h\leq i_h, 1\leq l_h\leq k_h$ for $1\leq h\leq
r$. As a consequence, we obtain that

$$T_s^r(p,\widetilde{M})=\underbrace{T_p\widetilde{M}\otimes\cdots\otimes
T_p\widetilde{M}}\limits_{r}\otimes\underbrace{T_p^*\widetilde{M}\otimes\cdots\otimes
T_p^*\widetilde{M}}\limits_{s}.$$

Since ${\rm dim}T_p\widetilde{M}={\rm
dim}T_p^*\widetilde{M}=\widehat{s}(p)+
\sum\limits_{i=1}^{s(p)}(n_i-\widehat{s}(p))$ by Theorems $3.2$ and
$3.3$, we also know that

$${\rm dim}T_s^r(p,\widetilde{M})=(\widehat{s}(p)+ \sum\limits_{i=1}^{s(p)}(n_i-\widehat{s}(p)))^{r+s}.$$

\vskip 4mm

\no{\bf Definition $3.7$} \ {\it Let
$T_s^r(\widetilde{M})=\bigcup\limits_{p\in\widetilde{M}}T_s^r(p,\widetilde{M})$
for a smoothly combinatorial manifold
$\widetilde{M}=\widetilde{M}(n_1,n_2,\cdots,n_m)$. A tensor filed
of type $(r,s)$ on $\widetilde{M}(n_1,n_2,\cdots,n_m)$ is a
mapping $\tau: \widetilde{M}(n_1,n_2,\cdots,n_m)\rightarrow
T_s^r(\widetilde{M})$ such that $\tau(p)\in
T_s^r(p,\widetilde{M})$ for $\forall
p\in\widetilde{M}(n_1,n_2,\cdots,n_m)$.

A $k$-form on $\widetilde{M}(n_1,n_2,\cdots,n_m)$ is a tensor field
$\omega\in T_0^r(\widetilde{M})$. Denoted all $k$-form of
$\widetilde{M}(n_1,n_2,\cdots,n_m)$ by $\Lambda^k(\widetilde{M})$
and
$\Lambda(\widetilde{M})=\bigoplus\limits_{k=0}^{\widehat{s}(p)-s(p)\widehat{s}(p)
+\sum_{i=1}^{s(p)}n_i}\Lambda^k(\widetilde{M})$,
$\mathscr{X}(\widetilde{M})=\bigcup\limits_{p\in\widetilde{M}}
\mathscr{X}_p$.}

\vskip 3mm

Similar to the classical differential geometry, we can also define
operations $\varphi\wedge\psi$ for $\forall \varphi, \psi\in
T_s^r(\widetilde{M})$, $[X,Y]$ for $\forall
X,Y\in\mathscr{X}(\widetilde{M})$ and obtain a {\it Lie algebra}
under the commutator. For the exterior differentiations on
combinatorial manifolds, we find results following.

\vskip 4mm

\no{\bf Theorem $3.5$}\ {\it Let $\widetilde{M}$ be a smoothly
combinatorial manifold. Then there is a unique exterior
differentiation
$\widetilde{d}:\Lambda(\widetilde{M})\rightarrow\Lambda(\widetilde{M})$
such that for any integer $k\geq 1$,
$\widetilde{d}(\Lambda^k)\subset\Lambda^{k+1}(\widetilde{M})$ with
conditions following hold.}

($1$) {\it $\widetilde{d}$ is linear, i.e., for $\forall \varphi,
\psi\in\Lambda(\widetilde{M})$, $\lambda\in{\bf R}$,}

$$\widetilde{d}(\varphi+\lambda\psi)= \widetilde{d}\varphi\wedge\psi+\lambda \widetilde{d}\psi$$

\no{\it and for $\varphi\in\Lambda^k(\widetilde{M}),
\psi\in\Lambda(\widetilde{M})$,}

$$\widetilde{d}(\varphi\wedge\psi)= \widetilde{d}\varphi+(-1)^k\varphi\wedge \widetilde{d}\psi.$$

($2$) {\it For $f\in\Lambda^0(\widetilde{M})$, $\widetilde{d}f$ is
the differentiation of $f$.}

($3$) {\it $\widetilde{d}^2=\widetilde{d}\cdot \widetilde{d}=0$.}

($4$) {\it $\widetilde{d}$ is a local operator, i.e., if $U\subset
V\subset\widetilde{M}$ are open sets and $\alpha\in\Lambda^k(V)$,
then $\widetilde{d}(\alpha|_U)= (\widetilde{d}\alpha)|_U$.}

\vskip 3mm

{\it Proof} \ Let $(U; [\varphi])$, where
$[\varphi]:p\rightarrow\bigcup\limits_{i=1}^{s(p)}$
$[\varphi](p)=[\varphi(p)]$ be a local chart for a point
$p\in\widetilde{M}$ and $\alpha=\alpha_{(\mu_1\nu_1)\cdots
(\mu_k\psi_k)}dx^{\mu_1\nu_1}\wedge\cdots\wedge dx^{\mu_k\nu_k}$
with $1\leq \nu_j\leq n_{\mu_i}$ for $1\leq\mu_i\leq s(p)$, $1\leq
i\leq k$. We first establish the uniqueness. If $k=0$, the local
formula $\widetilde{d}\alpha =\frac{\partial\alpha}{\partial
x^{\mu\nu}}dx^{\mu\nu}$ applied to the coordinates $x^{\mu\nu}$ with
$1\leq \nu_j\leq n_{\mu_i}$ for $1\leq\mu_i\leq s(p)$, $1\leq i\leq
k$ shows that the differential of $x^{\mu\nu}$ is $1$-form
$dx^{\mu\nu}$. From $(3)$, $\widetilde{d}(x^{\mu\nu})=0$, which
combining with $(1)$ shows that
$\widetilde{d}(dx^{\mu_1\nu_1}\wedge\cdots\wedge
d^{x^{\mu_k\nu_k}})=0$. This, again by $(1)$,

$$\widetilde{d}\alpha=\frac{\partial\alpha_{(\mu_1\nu_1)\cdots
(\mu_k\psi_k)}}{\partial x^{\mu\nu}}dx^{\mu\nu} \wedge
dx^{\mu_1\nu_1}\wedge\cdots\wedge dx^{\mu_k\nu_k}. \ \ \ (3.3)$$

\no and $\widetilde{d}$ is uniquely determined on $U$ by properties
$(1)-(3)$ and by $(4)$ on any open subset of $\widetilde{M}$.

For existence, define on every local chart $(U; [\varphi])$ the
operator $\widetilde{d}$ by $(3.3)$. Then ($2$) is trivially
verified as is ${\bf R}$-linearity. If
$\beta=\beta_{(\sigma_1\varsigma_1)\cdots
(\sigma_l\varsigma_l)}dx^{\sigma_1\varsigma_1}\wedge\cdots\wedge
dx^{\sigma_l\varsigma_l}\in\Lambda^l(U)$, then

\begin{eqnarray*}
\widetilde{d}(\alpha\wedge\beta)&=&
\widetilde{d}(\alpha_{(\mu_1\nu_1)\cdots
(\mu_k\psi_k)}\beta_{(\sigma_1\varsigma_1)\cdots
(\sigma_l\varsigma_l)}dx^{\mu_1\nu_1}\wedge\cdots\wedge
d^{x^{\mu_k\nu_k}}\wedge dx^{\sigma_1\varsigma_1}\wedge\cdots\wedge
dx^{\sigma_l\varsigma_l})\\
&=& (\frac{\partial\alpha_{(\mu_1\nu_1)\cdots
(\mu_k\psi_k)}}{\partial
x^{\mu\nu}}\beta_{(\sigma_1\varsigma_1)\cdots
(\sigma_l\varsigma_l)}+\alpha_{(\mu_1\nu_1)\cdots (\mu_k\psi_k)}\\
&\times& \frac{\partial\beta_{(\sigma_1\varsigma_1)\cdots
(\sigma_l\varsigma_l)}}{\partial
x^{\mu\nu}})dx^{\mu_1\nu_1}\wedge\cdots\wedge
d^{x^{\mu_k\nu_k}}\wedge dx^{\sigma_1\varsigma_1}\wedge\cdots\wedge
dx^{\sigma_l\varsigma_l}\\
&=& \frac{\partial\alpha_{(\mu_1\nu_1)\cdots
(\mu_k\psi_k)}}{\partial
x^{\mu\nu}}dx^{\mu_1\nu_1}\wedge\cdots\wedge
d^{x^{\mu_k\nu_k}}\wedge\beta_{(\sigma_1\varsigma_1)\cdots
(\sigma_l\varsigma_l)}dx^{\sigma_1\varsigma_1}\wedge\cdots\wedge
dx^{\sigma_l\varsigma_l}\\
&+& (-1)^k\alpha_{(\mu_1\nu_1)\cdots
(\mu_k\psi_k)}dx^{\mu_1\nu_1}\cdots\wedge
d^{x^{\mu_k\nu_k}}\wedge\frac{\partial\beta_{(\sigma_1\varsigma_1)\cdots
(\sigma_l\varsigma_l)}}{\partial
x^{\mu\nu}})dx^{\sigma_1\varsigma_1}\cdots\wedge
dx^{\sigma_l\varsigma_l}\\
&=&
\widetilde{d}\alpha\wedge\beta+(-1)^k\alpha\wedge\widetilde{d}\beta
\end{eqnarray*}

\no and ($1$) is verified. For ($3$), symmetry of the second partial
derivatives shows that

$$\widetilde{d}(\widetilde{d}\alpha)=\frac{\partial^2\alpha_{(\mu_1\nu_1)\cdots
(\mu_k\psi_k)}}{\partial x^{\mu\nu}\partial
x^{\sigma\varsigma}}dx^{\mu_1\nu_1}\wedge\cdots\wedge
d^{x^{\mu_k\nu_k}}\wedge dx^{\sigma_1\varsigma_1}\wedge\cdots\wedge
dx^{\sigma_l\varsigma_l}) = 0.$$

\no Thus, in every local chart $(U;[\varphi])$, ($3.3$) defines an
operator $\widetilde{d}$ satisfying ($1$)-($3$). It remains to be
shown that $\widetilde{d}$ really defines an operator
$\widetilde{d}$ on any open set and ($4$) holds. To do so, it
suffices to show that this definition is chart independent. Let
$\widetilde{d}'$ be the operator given by $(3.3)$ on a local chart
$(U';[\varphi'])$, where $U\bigcap U'\not=\emptyset$. Since
$\widetilde{d}'$ also satisfies $(1)-(3)$ and the local uniqueness
has already been proved, $\widetilde{d}'\alpha =\widetilde{d}\alpha$
on $U\bigcap U'$. Whence, ($4$) thus follows. \ \ \ $\natural$

\vskip 4mm

\no{\bf Corollary $3.2$} \ {\it Let
$\widetilde{M}=\widetilde{M}(n_1,n_2,\cdots,n_m)$ be a smoothly
combinatorial manifold and
$d_M:\Lambda^k(M)\rightarrow\Lambda^{k+1}(M)$ the unique exterior
differentiation on $M$ with conditions following hold for $M\in
V(G^l[\widetilde{M}(n_1,n_2,\cdots,n_m)])$ where, $1\leq
l\leq\min\{n_1,n_2,$ $\cdots,n_m\}$.}

($1$) {\it $d_M$ is linear, i.e., for $\forall \varphi,
\psi\in\Lambda(M)$, $\lambda\in{\bf R}$,}

$$d_M(\varphi+\lambda\psi)= d_M\varphi+\lambda d_M\psi.$$

($2$) {\it For $\varphi\in\Lambda^r(M), \psi\in\Lambda(M)$,}

$$d_M(\varphi\wedge\psi)= d_M\varphi+(-1)^r\varphi\wedge d_M\psi.$$

($3$) {\it For $f\in\Lambda^0(M)$, $d_Mf$ is the differentiation
of $f$.}

($4$) {\it $d_M^2=d_M\cdot d_M=0$.}

\no{\it Then }

 $$\widetilde{d}|_M= d_M.$$

 \vskip 3mm

 {\it Proof} \ By Theorem $2.4.5$ in $[1]$, $d_M$ exists uniquely for any
 smoothly manifold $M$. Now since $\widetilde{d}$ is a local
 operator on $\widetilde{M}$, i.e., for any open subset
 $U_{\mu}\subset\widetilde{M}$, $\widetilde{d}(\alpha|_{U_{\mu}})
 =(\widetilde{d}\alpha)|_{U_{\mu}}$ and there is
 an index set $J$ such that $M=\bigcup\limits_{\mu\in
 J}U_{\mu}$, we finally get that

 $$\widetilde{d}|_M= d_M$$

 \no by the uniqueness of $\widetilde{d}$ and $d_M$. \ \ \ $\natural$

\vskip 4mm

\no{\bf Theorem $3.6$} \ {\it Let
$\omega\in\Lambda^1(\widetilde{M})$. Then for $\forall
X,Y\in\mathscr{X}(\widetilde{M})$,}

$$\widetilde{d}\omega(X,Y)=X(\omega(Y))-Y(\omega(X))-\omega([X,Y]).$$

\vskip 3mm

{\it Proof} \ Denote by $\alpha(X,Y)$ the right hand side of the
formula. We know that
$\alpha:\widetilde{M}\times\widetilde{M}\rightarrow
C^{\infty}(\widetilde{M})$. It can be checked immediately that
$\alpha$ is bilinear and for $\forall
X,Y\in\mathscr{X}(\widetilde{M})$, $f\in C^{\infty}(\widetilde{M})$,

\begin{eqnarray*}
\alpha(fX,Y) &=& fX(\omega(Y))-Y(\omega(fX))-\omega([fX,Y])\\
&=& fX(\omega(Y))-Y(f\omega(X))-\omega(f[X,Y]-Y(f)X)\\
&=& f\alpha(X,Y)
\end{eqnarray*}

\no and

$$\alpha(X,fY)=-\alpha(fY,X)=-f\alpha(Y,X)=f\alpha(X,Y)$$

\no by definition. Accordingly, $\alpha$ is a differential $2$-form.
We only need to prove that for a local chart $(U,[\varphi])$,

$$\alpha|_U = \widetilde{d}\omega|_U.$$

\no In fact, assume $\omega|_U=\omega_{\mu\nu}dx^{\mu\nu}$. Then

\begin{eqnarray*}
(\widetilde{d}\omega)|_U=\widetilde{d}(\omega|_U) &=&
\frac{\partial\omega_{\mu\nu}}{\partial x^{\sigma\varsigma}} dx^{\sigma\varsigma}\wedge dx^{\mu\nu}\\
&=& \frac{1}{2}(\frac{\partial\omega_{\mu\nu}}{\partial
x^{\sigma\varsigma}}-\frac{\partial\omega_{\varsigma\tau}}{\partial
x^{\mu\nu}})dx^{\sigma\varsigma}\wedge dx^{\mu\nu}.
\end{eqnarray*}

\no On the other hand,
$\alpha|_U=\frac{1}{2}\alpha(\frac{\partial}{\partial
x^{\mu\nu}},\frac{\partial}{\partial
x^{\sigma\varsigma}})dx^{\sigma\varsigma}\wedge dx^{\mu\nu}$, where

\begin{eqnarray*}
\alpha(\frac{\partial}{\partial x^{\mu\nu}},\frac{\partial}{\partial
x^{\sigma\varsigma}}) &=& \frac{\partial}{\partial
x^{\sigma\varsigma}}(\omega(\frac{\partial}{\partial
x^{\mu\nu}}))-\frac{\partial}{\partial
x^{\mu\nu}}(\omega(\frac{\partial}{\partial x^{\sigma\varsigma}}))\\
& \ & -\omega([\frac{\partial}{\partial
x^{\mu\nu}}-\frac{\partial}{\partial x^{\sigma\varsigma}}])\\
&=& \frac{\partial\omega_{\mu\nu}}{\partial
x^{\sigma\varsigma}}-\frac{\partial\omega_{\sigma\varsigma}}{\partial
x^{\mu\nu}}.
\end{eqnarray*}

\no Therefore, $\widetilde{d}\omega|_U = \alpha|_U$. \ \ \
$\natural$

\vskip 4mm

\no{\bf $3.3$ Connections on tensors}

\vskip 3mm

\no We introduce connections on tensors of smoothly combinatorial
manifolds by the next definition.

\vskip 4mm

\no{\bf Definition $3.8$} \ {\it Let $\widetilde{M}$ be a smoothly
combinatorial manifold. A connection on tensors of $\widetilde{M}$
is a mapping $\widetilde{D}: \mathscr{X}(\widetilde{M})\times
T_s^r{\widetilde{M}}\rightarrow T_s^r{\widetilde{M}}$ with
$\widetilde{D}_X\tau=\widetilde{D}(X,\tau)$ such that for $\forall
X,Y\in\mathscr{X}{\widetilde{M}}$, $\tau, \pi\in
T_s^r(\widetilde{M})$,$\lambda\in{\bf R}$ and $f\in
C^{\infty}(\widetilde{M})$,}

($1$)
$\widetilde{D}_{X+fY}\tau=\widetilde{D}_X\tau+f\widetilde{D}_Y\tau$;
and $\widetilde{D}_X(\tau+\lambda
\pi)=\widetilde{D}_X\tau+\lambda\widetilde{D}_X\pi$;

($2$) $\widetilde{D}_X(\tau\otimes\pi)=\widetilde{D}_X\tau\otimes\pi
+\sigma\otimes\widetilde{D}_X\pi$;

($3$) {\it for any contraction $C$ on $T_s^r(\widetilde{M})$,}

$$\widetilde{D}_X(C(\tau))=C(\widetilde{D}_X\tau).$$

\vskip 3mm

We get results following for these connections on tensors of
smoothly combinatorial manifolds.

\vskip 4mm

\no{\bf Theorem $3.7$} \ {\it Let $\widetilde{M}$ be a smoothly
combinatorial manifold. Then there exists a connection
$\widetilde{D}$ locally on $\widetilde{M}$ with a form}

$$(\widetilde{D}_X\tau)|_U=X^{\sigma\varsigma}
\tau_{(\kappa_1\lambda_1)(\kappa_2\lambda_2)\cdots(\kappa_s\lambda_s),(\mu\nu)}
^{(\mu_1\nu_1)(\mu_2\nu_2)\cdots(\mu_r\nu_r)}\frac{\partial}{\partial
x^{\mu_1\nu_1}}\otimes\cdots\otimes\frac{\partial}{\partial
x^{\mu_r\nu_r}}\otimes dx^{\kappa_1\lambda_1}\otimes \cdots\otimes
dx^{\kappa_s\lambda_s}$$

\no{\it for $\forall Y\in\mathscr{X}(\widetilde{M})$ and $\tau\in
T_s^r(\widetilde{M})$, where}

\begin{eqnarray*}
\tau_{(\kappa_1\lambda_1)(\kappa_2\lambda_2)\cdots(\kappa_s\lambda_s),(\mu\nu)}
^{(\mu_1\nu_1)(\mu_2\nu_2)\cdots(\mu_r\nu_r)}&=&
\frac{\partial\tau_{(\kappa_1\lambda_1)(\kappa_2\lambda_2)\cdots
(\kappa_s\lambda_s)}^{(\mu_1\nu_1)(\mu_2\nu_2)\cdots(\mu_r\nu_r)}}{\partial
x^{\mu\nu}}\\
&+&
\sum\limits_{a=1}^r\tau_{(\kappa_1\lambda_1)(\kappa_2\lambda_2)\cdots(\kappa_s\lambda_s)}
^{(\mu_1\nu_1)\cdots(\mu_{a-1}\nu_{a-1})(\sigma\varsigma)(\mu_{a+1}\nu_{a+1})\cdots(\mu_r\nu_r)}
\Gamma_{(\sigma\varsigma)(\mu\nu)}^{\mu_a\nu_a}\\
&-& \sum\limits_{b=1}^s
\tau_{(\kappa_1\lambda_1)\cdots(\kappa_{b-1}\lambda_{b-1})(\mu\nu)(\sigma_{b+1}\varsigma_{b+1})
\cdots(\kappa_s\lambda_s)}
^{(\mu_1\nu_1)(\mu_2\nu_2)\cdots(\mu_r\nu_r)}\Gamma_{(\sigma_b\varsigma_b)(\mu\nu)}^{\sigma\varsigma}
\end{eqnarray*}

\no{\it and $\Gamma_{(\sigma\varsigma)(\mu\nu)}^{\kappa\lambda}$ is
a function determined by}

$$\widetilde{D}_{\frac{\partial}{\partial x^{\mu\nu}}}\frac{\partial}{\partial x^{\sigma\varsigma}}
=\Gamma_{(\sigma\varsigma)(\mu\nu)}^{\kappa\lambda}\frac{\partial}{\partial
x^{\sigma\varsigma}}$$

\no{\it on $(U_p;[\varphi_p])=(U_p;x^{\mu\nu})$ of a point
$p\in\widetilde{M}$, also called the coefficient on a connection.}

\vskip 3mm

{\it Proof} \ We first prove that any connection $\widetilde{D}$ on
smoothly combinatorial manifolds $\widetilde{M}$ is local by
definition, namely for $X_1,X_2\in\mathscr{X}(\widetilde{M})$ and
$\tau_1, \tau_2\in T_s^r(\widetilde{M})$, if $X_1|_U=X_2|_U$ and
$\tau_1|_U=\tau_2|_U$, then
$(\widetilde{D}_{X_1}\tau_1)_U=(\widetilde{D}_{X_2}\tau_2)_U$. For
this objective, we need to prove that
$(\widetilde{D}_{X_1}\tau_1)_U=(\widetilde{D}_{X_1}\tau_2)_U$ and
$(\widetilde{D}_{X_1}\tau_1)_U=(\widetilde{D}_{X_2}\tau_1)_U$. Since
their proofs are similar, we check the first only.

In fact, if $\tau=0$, then $\tau=\tau-\tau$. By the definition of
connection,

$$\widetilde{D}_X\tau=\widetilde{D}_X(\tau-\tau)
=\widetilde{D}_X\tau-\widetilde{D}_X\tau=0.$$

Now let $p\in U$. Then there is a neighborhood $V_p$ of $p$ such
that $\overline{V}$ is compact and $\overline{V}\subset U$. By a
result in topology, i.e., {\it for two open sets $V_p,U$ of ${\bf
R}^{\widehat{s}(p)-s(p)\widehat{s}(p)+n_1+\cdots+n_{s(p)}}$ with
compact $\overline{V_p}$ and $\overline{V_p}\subset U$, there exists
a function $f\in C^{\infty}({\bf
R}^{\widehat{s}(p)-s(p)\widehat{s}(p)+n_1+\cdots+n_{s(p)}})$ such
that $0\leq f\leq 1$ and $f|_{V_p}\equiv 1$, $f|_{{\bf
R}^{\widehat{s}(p)-s(p)\widehat{s}(p)+n_1+\cdots+n_{s(p)}}\setminus
U}\equiv 0$,} we find that $f\cdot(\tau_2-\tau_1)=0$. Whence, we
know that

$$0=\widetilde{D}_{X_1}((f\cdot(\tau_2-\tau_1)))
=X_1(f)(\tau_2-\tau_1)+f(\widetilde{D}_{X_1}\tau_2-\widetilde{D}_{X_1}\tau_1).$$

\no As a consequence, we get that
$(\widetilde{D}_{X_1}\tau_1)_V=(\widetilde{D}_{X_1}\tau_2)_V$,
particularly,
$(\widetilde{D}_{X_1}\tau_1)_p=(\widetilde{D}_{X_1}\tau_2)_p$. For
the arbitrary choice of $p$, we get that
$(\widetilde{D}_{X_1}\tau_1)_U=(\widetilde{D}_{X_1}\tau_2)_U$
finally.

The local property of $\widetilde{D}$ enables us to find an induced
connection $\widetilde{D}^U:\mathscr{X}(U)\times T_s^r(U)\rightarrow
T_s^r(U)$ such that
$\widetilde{D}_{X|_U}^U(\tau|_U)=(\widetilde{D}_X\tau)|_U$ for
$\forall X\in\mathscr{X}(\widetilde{M})$ and $\tau\in
T_s^r\widetilde{M}$. Now for $\forall
X_1,X_2\in\mathscr{X}(\widetilde{M})$, $\forall\tau_1, \tau_2\in
T_s^r(\widetilde{M})$ with $X_1|_{V_p}=X_2|_{V_p}$ and
$\tau_1|_{V_p}=\tau_2|_{V_p}$, define a mapping
$\widetilde{D}^{U}:\mathscr{X}(U)\times T_s^r(U)\rightarrow
T_s^r(U)$ by

$$(\widetilde{D}_{X_1}\tau_1)|_{V_p}=(\widetilde{D}_{X_1}\tau_2)|_{V_p}$$

\no for any point $p\in U$. Then since $\widetilde{D}$ is a
connection on $\widetilde{M}$, it can be checked easily that
$\widetilde{D}^U$ satisfies all conditions in Definition $3.8$.
Whence, $\widetilde{D}^U$ is indeed a connection on $U$.

Now we calculate the local form on a chart $(U_p,[\varphi_p])$ of
$p$. Since

$$\widetilde{D}_{\frac{\partial}{\partial x^{\mu\nu}}}
=\Gamma_{(\sigma\varsigma)(\mu\nu)}^{\kappa\lambda}\frac{\partial}{\partial
x^{\sigma\varsigma}},$$

\no it can find immediately that

$$
\widetilde{D}_{\frac{\partial}{\partial x^{\mu\nu}}}
dx^{\kappa\lambda}
=-\Gamma_{(\sigma\varsigma)(\mu\nu)}^{\kappa\lambda}
dx^{\sigma\varsigma}
$$

\no by Definition $3.8$. Therefore, we find that

$$(\widetilde{D}_X\tau)|_U=X^{\sigma\varsigma}
\tau_{(\kappa_1\lambda_1)(\kappa_2\lambda_2)\cdots(\kappa_s\lambda_s),(\mu\nu)}
^{(\mu_1\nu_1)(\mu_2\nu_2)\cdots(\mu_r\nu_r)}\frac{\partial}{\partial
x^{\mu_1\nu_1}}\otimes\cdots\otimes\frac{\partial}{\partial
x^{\mu_r\nu_r}}\otimes dx^{\kappa_1\lambda_1}\otimes \cdots\otimes
dx^{\kappa_s\lambda_s}$$

\no with

\begin{eqnarray*}
\tau_{(\kappa_1\lambda_1)(\kappa_2\lambda_2)\cdots(\kappa_s\lambda_s),(\mu\nu)}
^{(\mu_1\nu_1)(\mu_2\nu_2)\cdots(\mu_r\nu_r)}&=&
\frac{\partial\tau_{(\kappa_1\lambda_1)(\kappa_2\lambda_2)\cdots
(\kappa_s\lambda_s)}^{(\mu_1\nu_1)(\mu_2\nu_2)\cdots(\mu_r\nu_r)}}{\partial
x^{\mu\nu}}\\
&+&
\sum\limits_{a=1}^r\tau_{(\kappa_1\lambda_1)(\kappa_2\lambda_2)\cdots(\kappa_s\lambda_s)}
^{(\mu_1\nu_1)\cdots(\mu_{a-1}\nu_{a-1})(\sigma\varsigma)(\mu_{a+1}\nu_{a+1})\cdots(\mu_r\nu_r)}
\Gamma_{(\sigma\varsigma)(\mu\nu)}^{\mu_a\nu_a}\\
&-& \sum\limits_{b=1}^s
\tau_{(\kappa_1\lambda_1)\cdots(\kappa_{b-1}\lambda_{b-1})(\mu\nu)(\sigma_{b+1}\varsigma_{b+1})
\cdots(\kappa_s\lambda_s)}
^{(\mu_1\nu_1)(\mu_2\nu_2)\cdots(\mu_r\nu_r)}\Gamma_{(\sigma_b\varsigma_b)
(\mu\nu)}^{\sigma\varsigma}.
\end{eqnarray*}

\no This completes the proof. \ \ \ $\natural$

\vskip 4mm

\no{\bf Theorem $3,8$} \ {\it Let $\widetilde{M}$ be a smoothly
combinatorial manifold with a connection $\widetilde{D}$. Then for
$\forall X,Y\in\mathscr{X}(\widetilde{M})$,}

$$\widetilde{T}(X,Y)=\widetilde{D}_XY-\widetilde{D}_YX-[X,Y]$$

\no{\it is a tensor of type $(1,2)$ on $\widetilde{M}$.}

\vskip 3mm

{\it Proof} \ By definition, it is clear that
$\widetilde{T}:\mathscr{X}(\widetilde{M})\times\mathscr{X}(\widetilde{M})
\rightarrow\mathscr{X}(\widetilde{M})$ is antisymmetrical and
bilinear. We only need to check it is also linear on each element in
$C^{\infty}(\widetilde{M})$ for variables $X$ or $Y$. In fact, for
$\forall f\in C^{\infty}(\widetilde{M})$,

\begin{eqnarray*}
\widetilde{T}(fX,Y)&=&
\widetilde{D}_{fX}Y-\widetilde{D}_Y(fX)-[fX,Y]\\
&=& f\widetilde{D}_{X}Y-(Y(f)X+f\widetilde{D}_YX)\\
&-& (f[X,Y]-Y(f)X)= f\widetilde{T}(X,Y).
\end{eqnarray*}

\no and

$$\widetilde{T}(X,fY)= -\widetilde{T}(fY,X)=-f\widetilde{T}(Y,X)
=f\widetilde{T}(X,Y). \ \ \natural$$

Notice that

\begin{eqnarray*}
T(\frac{\partial}{\partial x^{\mu\nu}}, \frac{\partial}{\partial
x^{\sigma\varsigma}}) &=& \widetilde{D}_{\frac{\partial}{\partial
x^{\mu\nu}}}\frac{\partial}{\partial x^{\sigma\varsigma}}
-\widetilde{D}_{\frac{\partial}{\partial x^{\sigma\varsigma}}}
\frac{\partial}{\partial x^{\mu\nu}}\\
&=& (\Gamma_{(\mu\nu)(\sigma\varsigma)}^{\kappa\lambda}
-\Gamma_{(\sigma\varsigma)(\mu\nu)}^{\kappa\lambda})\frac{\partial}{\partial
x^{\kappa\lambda}}
\end{eqnarray*}

\no under a local chart $(U_p;[\varphi_p])$ of a point
$p\in\widetilde{M}$. If $T(\frac{\partial}{\partial x^{\mu\nu}},
\frac{\partial}{\partial x^{\sigma\varsigma}})\equiv 0$, we call $T$
{\it torsion-free}. This enables us getting the next useful result.

\vskip 4mm

\no{\bf Theorem $3.9$} \ {\it A connection $\widetilde{D}$ on
tensors of a smoothly combinatorial manifold $\widetilde{M}$ is
torsion-free if and only if
$\Gamma_{(\mu\nu)(\sigma\varsigma)}^{\kappa\lambda}
=\Gamma_{(\sigma\varsigma)(\mu\nu)}^{\kappa\lambda}$.}

\vskip 3mm

Now we turn our attention to the case of $s=r=1$. Similarly, a {\it
combinatorially Riemannian geometry} is defined in the next
definition.

\vskip 4mm

\no{\bf Definition $3.9$} \ {\it Let $\widetilde{M}$ be a smoothly
combinatorial manifold and $g\in
A^2(\widetilde{M})=\bigcup\limits_{p\in\widetilde{M}}T_2^0(p,\widetilde{M})$.
If $g$ is symmetrical and positive, then $\widetilde{M}$ is called a
combinatorially Riemannian manifold, denoted by $(\widetilde{M},g)$.
In this case, if there is a connection $\widetilde{D}$ on
$(\widetilde{M},g)$ with equality following hold

$$Z(g(X,Y))=g(\widetilde{D}_Z,Y)+g(X,\widetilde{D}_ZY) \ \ \ (3.4)$$

\no then $\widetilde{M}$ is called a combinatorially Riemannian
geometry, denoted by $(\widetilde{M},g, \widetilde{D})$.}

\vskip 3mm

We get a result for connections on smoothly combinatorial manifolds
similar to that of Riemannian geometry.

\vskip 4mm

\no{\bf Theorem $3.10$} \ {\it Let $(\widetilde{M},g)$ be a
combinatorially Riemannian manifold. Then there exists a unique
connection $\widetilde{D}$ on $(\widetilde{M},g)$ such that
$(\widetilde{M},g, \widetilde{D})$ is a combinatorially Riemannian
geometry.}

\vskip 3mm

{\it Proof} \ By definition, we know that

$$\widetilde{D}_Zg(X,Y)=Z(g(X,Y))-g(\widetilde{D}_ZX,Y)-g(X,\widetilde{D}_ZY)$$

\no for a connection $\widetilde{D}$ on tensors of $\widetilde{M}$
and $\forall Z\in\mathscr{X}(\widetilde{M})$. Thereby, the equality
$(3.4)$ is equivalent to that of $\widetilde{D}_Zg=0$ for $\forall
Z\in\mathscr{X}(\widetilde{M})$, namely $\widetilde{D}$ is
torsion-free.

Not loss of generality, assume
$g=g_{(\mu\nu)(\sigma\varsigma)}dx^{\mu\nu}dx^{\sigma\varsigma}$ in
a local chart $(U_p;[\varphi_p])$ of a point $p$, where
$g_{(\mu\nu)(\sigma\varsigma)}=g(\frac{\partial}{\partial
x^{\mu\nu}},\frac{\partial}{\partial x^{\sigma\varsigma}})$. Then we
find that

$$\widetilde{D}g=(\frac{\partial g_{(\mu\nu)(\sigma\varsigma)}}
{\partial x^{\kappa\lambda}}-g_{(\zeta\eta)(\sigma\varsigma)}
\Gamma_{(\mu\nu)(\sigma\varsigma)}^{\zeta\eta}
-g_{(\mu\nu)(\zeta\eta)}\Gamma_{(\sigma\varsigma)(\kappa\lambda)}^{\zeta\eta})
dx^{\mu\nu}\otimes dx^{\sigma\varsigma}\otimes dx^{\kappa\lambda}.$$

\no Therefore, we get that

$$\frac{\partial g_{(\mu\nu)(\sigma\varsigma)}}
{\partial x^{\kappa\lambda}}=g_{(\zeta\eta)(\sigma\varsigma)}
\Gamma_{(\mu\nu)(\sigma\varsigma)}^{\zeta\eta}
+g_{(\mu\nu)(\zeta\eta)}\Gamma_{(\sigma\varsigma)(\kappa\lambda)}^{\zeta\eta}
\ \ \ (3.5)$$

\no if $\widetilde{D}_Zg=0$ for $\forall
Z\in\mathscr{X}(\widetilde{M})$. The formula $(3.5)$ enables us to
get that

$$\Gamma_{(\mu\nu)(\sigma\varsigma)}^{\kappa\lambda}
=\frac{1}{2}g^{(\kappa\lambda)(\zeta\eta)}(\frac{\partial
g_{(\mu\nu)(\zeta\eta)}}{\partial
x^{\sigma\varsigma}}+\frac{\partial
g_{(\zeta\eta)(\sigma\varsigma)}}{\partial
x^{\mu\nu}}-\frac{\partial g_{(\mu\nu)(\sigma\varsigma)}}{\partial
x^{\zeta\eta}}),$$

\no where $g^{(\kappa\lambda)(\zeta\eta)}$ is an element in the
matrix inverse of $[g_{(\mu\nu)(\sigma\varsigma)}]$.

Now if there exists another torsion-free connection
$\widetilde{D}^*$ on $(\widetilde{M},g)$ with

$$\widetilde{D}^*_{\frac{\partial}{\partial x^{\mu\nu}}}
={\Gamma^*}_{(\sigma\varsigma)(\mu\nu)}^{\kappa\lambda}\frac{\partial}{\partial
x^{\kappa\lambda}},$$

\no then we must get that

$${\Gamma^*}_{(\mu\nu)(\sigma\varsigma)}^{\kappa\lambda}
=\frac{1}{2}g^{(\kappa\lambda)(\zeta\eta)}(\frac{\partial
g_{(\mu\nu)(\zeta\eta)}}{\partial
x^{\sigma\varsigma}}+\frac{\partial
g_{(\zeta\eta)(\sigma\varsigma)}}{\partial
x^{\mu\nu}}-\frac{\partial g_{(\mu\nu)(\sigma\varsigma)}}{\partial
x^{\zeta\eta}}).$$

\no Accordingly, $\widetilde{D}=\widetilde{D}^*$. Whence, there are
at most one torsion-free connection $\widetilde{D}$ on a
combinatorially Riemannian manifold $(\widetilde{M},g)$.

For the existence of torsion-free connection $\widetilde{D}$ on
$(\widetilde{M},g)$, let
$\Gamma_{(\mu\nu)(\sigma\varsigma)}^{\kappa\lambda}
=\Gamma_{(\sigma\varsigma)(\mu\nu)}^{\kappa\lambda}$ and define a
connection $\widetilde{D}$ on  $(\widetilde{M},g)$ such that

$$\widetilde{D}_{\frac{\partial}{\partial x^{\mu\nu}}}
=\Gamma_{(\sigma\varsigma)(\mu\nu)}^{\kappa\lambda}\frac{\partial}{\partial
x^{\kappa\lambda}},$$

\no then $\widetilde{D}$ is torsion-free by Theorem $3.9$. This
completes the proof. \ \ \ $\natural$

\vskip 4mm

\no{\bf Corollary $3.3$}([2]) \ {\it For a Riemannian manifold
$(M,g)$, there exists only one torsion-free connection $D$, i.e.,}

$$D_Zg(X,Y)=Z(g(X,Y))-g(D_ZX,Y)-g(X,D_ZY)\equiv 0$$

\no{\it for $\forall X,Y,Z\in\mathscr{X}(M)$.}

\vskip 4mm

\no{\bf $3.4$ Minkowski Norms}

\vskip 3mm

\no These {\it Minkowski norms} are the fundamental in Finsler
geometry. Certainly, they can be also generalized on smoothly
combinatorial manifolds.

\vskip 4mm

\no{\bf Definition $3.10$} \ {\it A Minkowski norm on a vector space
$V$ is a function $F:V\rightarrow {\bf R}$ such that}

$(1)$ \ {\it $F$ is smooth on $V\backslash\{0\}$ and $F(v)\geq 0$
for $\forall v\in V$;}

$(2)$ \ {\it $F$ is $1$-homogenous, i.e., $F(\lambda v)=\lambda
F(v)$ for $\forall\lambda> 0$;}

$(3)$ \ {\it for all $y\in V\backslash\{0\}$, the symmetric bilinear
form $g_y: V\times V\rightarrow{\bf R}$ with}

$$g_y(u,v)=\sum\limits_{i,j}\frac{\partial^2 F(y)}{\partial y^i\partial y^j}$$

\no{\it is positive definite for $u,v\in V$.}

\vskip 3mm

Denoted by
$T\widetilde{M}=\bigcup\limits_{p\in\widetilde{M}}T_p\widetilde{M}$.
Similar to Finsler geometry, we introduce combinatorially Finsler
geometries on a Minkowski norm defined on $T\widetilde{M}$.

\vskip 4mm

\no{\bf Definition $3.11$} \ {\it A combinatorially Finsler geometry
is a smoothly combinatorial manifold $\widetilde{M}$ endowed with a
Minkowski norm $\widetilde{F}$ on $T\widetilde{M}$, denoted by
$(\widetilde{M}; \widetilde{F})$.}

\vskip 3mm

Then we get the following result.

\vskip 4mm

\no{\bf Theorem $3.11$} \ {\it There are combinatorially Finsler
geometries.}

\vskip 3mm

{\it Proof} \ Let $\widetilde{M}(n_1,n_2,\cdots,n_m)$ be a smoothly
combinatorial manifold. We construct Minkowski norms on
$T\widetilde{M}(n_1,n_2,\cdots,n_m)$. Let ${\bf
R}^{n_1+n_2+\cdots+n_m}$ be an eucildean space. Then there exists a
Minkowski norm $F(\overline{x})=|\overline{x}|$ in ${\bf
R}^{n_1+n_2+\cdots+n_m}$ at least, in here $|\overline{x}|$ denotes
the euclidean norm on ${\bf R}^{n_1+n_2+\cdots+n_m}$. According to
Theorem $3.2$, $T_p\widetilde{M}(n_1,n_2,\cdots,n_m)$ is
homeomorphic to ${\bf
R}^{\widehat{s}(p)-s(p)\widehat{s}(p)+n_{i_1}+\cdots+n_{i_{s(p)}}}$.
Whence there are Minkowski norms on
$T_p\widetilde{M}(n_1,n_2,\cdots,n_m)$ for $p\in U_p$, where
$(U_p;[\varphi_p])$ is a local chart.

Notice that the number of manifolds are finite in a smoothly
combinatorial manifold $\widetilde{M}(n_1,n_2,\cdots,n_m)$ and each
manifold has a finite cover
$\{(U_{\alpha};\varphi_{\alpha})|\alpha\in I\}$, where $I$ is a
finite index set. We know that there is a finite cover

$$\bigcup\limits_{M\in
V(G[\widetilde{M}(n_1,n_2,\cdots,n_m)])}\{(U_{M\alpha};\varphi_{M\alpha})|\alpha\in
I_M\}.$$

\no By the decomposition theorem for unit, we know that there are
smooth functions $h_{M\alpha}, \alpha\in I_M$ such that

$$\sum\limits_{M\in V(G[\widetilde{M}(n_1,n_2,\cdots,n_m)])}
\sum\limits_{\alpha\in I_M}h_{M\alpha}=1 \ {\rm with} \ 0\leq h_{M\alpha}\leq 1.$$

Now we choose a Minkowski norm $\widetilde{F}^{M\alpha}$ on
$T_pM_{\alpha}$ for $\forall p\in U_{M\alpha}$. Define

\[
\widetilde{F}_{M\alpha}=\left\{\begin{array}{cc}
h^{M\alpha}\widetilde{F}^{M\alpha},& {\rm if}\quad p\in U_{M\alpha} ,\\
0,& {\rm if}\quad p\not\in U_{M\alpha}
\end{array}
\right.
\]

\no for $\forall p\in\widetilde{M}$. Now let

$$\widetilde{F}=\sum\limits_{M\in V(G[\widetilde{M}(n_1,n_2,\cdots,n_m)])}
\sum\limits_{\alpha\in I}\widetilde{F}_{M\alpha}.$$

\no Then $\widetilde{F}$ is a Minkowski norm on
$T\widetilde{M}(n_1,n_2,\cdots,n_m)$ since it can be checked
immediately that all conditions $(1)-(3)$ in Definition $3.10$ hold.
\ \ $\natural$

For the relation of combinatorially Finsler geometries with these
Smarandache geometries, we obtain the next consequence.

\vskip 4mm

\no{\bf Theorem $3.12$} \ {\it A combinatorially Finsler geometry
$(\widetilde{M}(n_1,n_2,\cdots,n_m); \widetilde{F})$ is a
Smarandache geometry if $m\geq 2$.}

\vskip 3mm

{\it Proof} \ Notice that if $m\geq 2$, then
$\widetilde{M}(n_1,n_2,\cdots,n_m)$ is combined by at least two
manifolds $M^{n_1}$ and $M^{n_2}$ with $n_1\not=n_2$. By definition,
we know that

$$M^{n_1}\setminus M^{n_2}\not=\emptyset \ {\rm and} \ M^{n_2}\setminus M^{n_1}\not=\emptyset.$$

\no Now the axiom {\it there is an integer $n$ such that there
exists a neighborhood homeomorphic to a open ball $B^n$ for any
point in this space} is Smarandachely denied, since for points in
$M^{n_1}\setminus M^{n_2}$, each has a neighborhood homeomorphic to
$B^{n_1}$, but each point in $M^{n_2}\setminus M^{n_1}$ has a
neighborhood homeomorphic to $B^{n_2}$. \ \ \ $\natural$

Theorems $3.11$ and $3.12$ imply inclusions in Smarandache
geometries for classical geometries in the following.

\vskip 4mm

\no{\bf Corollary $3.5$} \ {\it There are inclusions among
Smarandache geometries, Finsler geometry, Riemannian geometry and
Weyl geometry}:

\begin{eqnarray*}
& \ &\{Smarandache \ geometries\}\supset\{combinatorially \ Finsler
\ geometries\}\\
& \ &\supset\{Finsler \ geometry\} \ and \ \{combinatorially \
Riemannian \ geometries\}\\
& \ &\supset\{Riemannian \ geometry\}\supset\{Weyl \ geometry\}.
\end{eqnarray*}

\vskip 3mm

{\it Proof} \ Let $m=1$. Then a combinatorially Finsler geometry
$(\widetilde{M}(n_1,n_2,\cdots,n_m); \widetilde{F})$ is nothing but
just a Finsler geometry. Applying Theorems $3.11$ and $3.12$ to this
special case, we get these inclusions as expected. \ \ \ $\natural$

\vskip 4mm

\no{\bf Corollary $3.6$}\ {\it There are inclusions among
Smarandache geometries, combinatorially Riemannian geometries and
K\"{a}hler geometry}:

\begin{eqnarray*}
\{Smarandache \ geometries\}&\supset&\{combinatorially \ Riemannian \ geometries\}\\
&\supset&\{Riemannian \ geometry\}\\
&\supset& \{K\ddot{a}hler \ geometry\}.
\end{eqnarray*}

\vskip 3mm

{\it Proof} \ Let $m=1$ in a combinatorial manifold
$\widetilde{M}(n_1,n_2,\cdots,n_m)$ and applies Theorems $3.10$ and
$3.12$, we get inclusions

\begin{eqnarray*}
\{\rm Smarandache \ geometries\}&\supset&\{\rm combinatorially \ Riemannian \ geometries\}\\
&\supset& \{\rm Riemannian \ geometry\}.
\end{eqnarray*}

For the K\"{a}hler geometry, notice that any complex manifold
$M^{n}_c$ is equal to a smoothly real manifold $M^{2n}$ with a
natural base $\{\frac{\partial}{\partial x^i},
\frac{\partial}{\partial y^i}\}$ for $T_pM^n_c$ at each point $p\in
M^n_c$. Whence, we get

$$\{\rm Riemannian \ geometry\}\supset\{\rm K\ddot{a}hler \
geometry\}. \ \ \ \natural$$

\vskip 6mm

\no{\bf \S $4.$ Further Discussions}

\vskip 3mm

\no{\bf $4.1$ Embedding problems} \ Whitney had shown that {\it any
smooth manifold $M^d$ can be embedded as a closed submanifold of
${\bf R}^{2d+1}$} in $1936$ ([1]). The same embedding problem for
finitely combinatorial manifold in an euclidean space is also
interesting. Since $\widetilde{M}$ is finite, by applying Whitney
theorem, we know that there is an integer $n(\widetilde{M}),
n(\widetilde{M}) < +\infty$ such that $\widetilde{M}$ can be
embedded as a closed submanifold in ${\bf R}^{n(\widetilde{M})}$.
Then {\it what is the minimum dimension of euclidean spaces
embeddable a given finitely combinatorial manifold $\widetilde{M}$}?
{\it Wether can we determine it for some combinatorial manifolds
with a given graph structure, such as those of complete graphs
$K^n$, circuits $P^n$ or cubic graphs $Q^n$}?

\vskip 4mm

\no{\bf Conjecture $4.1$} \ {\it The minimum dimension of euclidean
spaces embeddable a finitely combinatorial manifold $\widetilde{M}$
is}

$$2\min\limits_{p\in\widetilde{M}}\{\widehat{s}(p)-s(p)\widehat{s}(p)+n_{i_1}+n_{i_2}+\cdots+n_{i_{s(p)}}\}+1.$$

\vskip 3mm

\no{\bf $4.2$ $D$-dimensional holes} \ For these closed
$2$-manifolds $S$, it is well-known that

\[
\chi (S)=\left\{\begin{array}{lr}
2-2p(S),& {\rm if} \ S \ {\rm is \ orientable} ,\\
2-q(S).& {\rm if} \ S {\rm is \ non-orientable.}
\end{array}
\right.
\]

\no with $p(S)$ or $q(S)$ the orientable genus or non-orientable
genus of $S$, namely $2$-dimensional holes adjacent to $S$. For
general case of $n$-manifolds $M$, we know that

$$\chi(M)=\sum\limits_{k=0}^{\infty}(-1)^k{\rm dim}H_k(M),$$

\no where ${\rm dim}H_k(M)$ is the rank of these $k$-dimensional
homolopy groups $H_k(M)$ in $M$, namely the number of
$k$-dimensional holes adjacent to the manifold $M$. By the
definition of combinatorial manifolds, some $k$-dimensional holes
adjacent to a combinatorial manifold are increased. Then {\it what
is the relation between the Euler-Poincare characteristic of a
combinatorial manifold $\widetilde{M}$ and the $i$-dimensional holes
adjacent to $\widetilde{M}$}? {\it Wether can we find a formula
likewise the Euler-Poincare formula}? Calculation shows that even
for the case of $n=2$, the situation is complex. For example, choose
$n$ different orientable $2$-manifolds $S_1,S_2,\cdots, S_n$ and let
them intersects one after another at $n$ different points in ${\bf
R}^3$. We get a combinatorial manifold $\widetilde{M}$. Calculation
shows that

$$\chi(\widetilde{M})=(\chi(S_1)+\chi(S_2)+\cdots+\chi(S_n))-n$$

\no by Theorem $2.9$. But it only increases one $2$-holes. {\it What
is the relation of $2$-dimensional holes adjacent to
$\widetilde{M}$}?

\vskip 3mm

\no{\bf $4.3$ Local properties} \ Although a finitely combinatorial
manifold $\widetilde{M}$ is not homogenous in general, namely the
dimension of local charts of two points in $\widetilde{M}$ maybe
different, we have still constructed global operators such as those
of exterior differentiation $\widetilde{d}$ and connection
$\widetilde{D}$ on $T_s^r\widetilde{M}$. A operator
$\widetilde{\mathfrak{O}}$ is said to be {\it local} on a subset
$W\subset T_s^r\widetilde{M}$ if for any local chart
$(U_p,[\varphi_p])$ of a point $p\in W$,

$$\widetilde{\mathfrak{O}}|_{U_p}(W)=\widetilde{\mathfrak{O}}(W)_{U_p}.$$

Of course, nearly all existent operators with local properties on
$T_s^r\widetilde{M}$ in Finsler or Riemannian geometries can be
reconstructed in these combinatorially Finsler or Riemannian
geometries and find the local forms similar to those in Finsler or
Riemannian geometries.

\vskip 3mm

\no{\bf $4.4$ Global properties} \ To find global properties on
manifolds is a central task in classical differential geometry. The
same is true for combinatorial manifolds. In classical geometry on
manifolds, some global results, such as those of {\it de Rham}
theorem and {\it Atiyah-Singer} index theorem,..., etc. are
well-known. Remember that the $p^{th}$ de Rham cohomology group on a
manifold $M$ and the {\it index} ${\rm Ind}\mathcal {D}$ of a {\it
Fredholm} operator $\mathcal {D}:H^k(M,E)\rightarrow L^2(M,F)$ are
defined to be a quotient space

$$H^p(M)=\frac{Ker(d:\Lambda^p(M)\rightarrow\Lambda^{p+1}(M))}
{Im(d:\Lambda^{p-1}(M)\rightarrow\Lambda^{p}(M))}.$$

\no and an integer

$${\rm Ind}\mathcal{D}={\rm dim}{Ker(\mathcal{D})}
-dim(\frac{L^2(M,F)}{{\rm Im}\mathcal{D}})$$

\no respectively. The de Rham theorem and the Atiyah-Singer index
theorem respectively conclude that\vskip 3mm

{\it for any manifold $M$, a mapping
$\varphi:\Lambda^p(M)\rightarrow Hom(\Pi_p(M),{\bf R})$ induces a
natural isomorphism $\varphi^*: H^p(M)\rightarrow H^n(M;{\bf R})$ of
cohomology groups, where $\Pi_p(M)$ is the free Abelian group
generated by the set of all $p$-simplexes in $M$}\vskip 2mm

\no and

$${\rm Ind}\mathcal{D}=Ind_T(\sigma(\mathcal {D})),$$

\no where $\sigma(\mathcal {D})):T^*M\rightarrow Hom(E,F)$
 and $Ind_T(\sigma(\mathcal {D}))$ is the topological index of
$\sigma(\mathcal D)$. Now the questions for these finitely
combinatorial manifolds are given in the following.

\vskip 3mm

($1$) {\it Is the de Rham theorem and Atiyah-Singer index theorem
still true for finitely combinatorial manifolds? If not, what is its
modified forms?}

($2$) {\it Check other global results for manifolds whether true or
get their new modified forms for finitely combinatorial manifolds.}

\vskip 10mm

\no{\bf References}\vskip 3mm

\re{[1]}R.Abraham and Marsden, {\it Foundation of Mechanics}(2nd
edition), Addison-Wesley, Reading, Mass, 1978.

\re{[2]}W.H.Chern and X.X.Li, {\it Introduction to Riemannian
Geometry}, Peking University Press, 2002.

\re{[3]}H.Iseri, {\it Smarandache manifolds}, American Research
Press, Rehoboth, NM,2002.

\re{[4]}L.Kuciuk and M.Antholy, An Introduction to Smarandache
Geometries, {\it Mathematics Magazine, Aurora, Canada},
Vol.12(2003).

\re{[5]}L.F.Mao, {\it Automorphism Groups of Maps, Surfaces and
Smarandache Geometries}, American Research Press, 2005.

\re{[6]}L.F.Mao, {\it On Automorphisms groups of Maps, Surfaces and
Smarandache geometries}, {\it Sientia Magna}, Vol.$1$(2005), No.$2$,
55-73.

\re{[7]}L.F.Mao, {\it Smarandache multi-space theory}, Hexis,
Phoenix, AZ£¬2006.

\re{[8]}L.F.Mao, On multi-metric spaces, {\it Scientia Magna},
Vol.2,No.1(2006), 87-94.

\re{[9]}L.F.Mao, On algebraic multi-group spaces, {\it Scientia
Magna}, Vol.2,No.1(2006), 64-70.

\re{[10]}L.F.Mao, On algebraic multi-ring spaces, {\it Scientia
Magna}, Vol.2,No.2(2006), 48-54.

\re{[11]}L.F.Mao, On algebraic multi-vector spaces, {\it Scientia
Magna}, Vol.2,No.2(2006), 1-6.

\re{[12]}L.F.Mao, Pseudo-Manifold Geometries with Applications,
e-print: {\it arXiv: math.} {\it GM/0610307}.

\re{[13]}L.F.Mao, A new view of combinatorial maps by Smarandache's
notion, {\it arXiv: math.GM/0506232}, also in {\it Selected Papers
on Mathematical Combinatorics}(I), World Academic Union, 2006.

\re{[14]}W.S.Massey, {\it Algebraic topology: an introduction},
Springer-Verlag,New York, etc.(1977).

\re{[15]}V.V.Nikulin and I.R.Shafarevlch, {\it Geometries and
Groups}, Springer-Verlag Berlin Heidelberg (1987)

\re{[16]}Joseph J.Rotman, {\it An introduction to algebraic
topology}, Springer-Verlag New York Inc. 1988.

\re{[17]}F.Smarandache, Mixed noneuclidean geometries, {\it eprint
arXiv: math/0010119}, 10/2000.

\re{[18]}J.Stillwell, {\it Classical topology and combinatorial
group theory}, Springer-Verlag New York Inc., (1980).

\end{document}